\documentclass[alpha-refs]{wiley-article}

\usepackage{color}
\usepackage{ulem}

\usepackage{siunitx}

\papertype{Original Article}
\paperfield{Journal Section}

\graphicspath{{figures/}}
\usepackage[pdfencoding=auto]{hyperref}
\usepackage{graphicx}%
\usepackage{multirow}%
\usepackage{amsmath,amssymb,amsfonts}%

\usepackage{amsthm}%
\usepackage{mathrsfs}%
\usepackage[title]{appendix}%
\usepackage{xcolor}%
\usepackage{textcomp}%
\usepackage{manyfoot}%
\usepackage{booktabs}%
\usepackage{algorithm}%
\usepackage{algorithmicx}
\usepackage{setspace}
\usepackage{algpseudocode}%
\usepackage{listings}%

\newtheorem{res}{Result}

\title{Estimating Counts Through an Average Rounded to the Nearest Non-negative Integer and its Theoretical \& Practical Effects}

\author[1\authfn{1}]{Roberto Rivera}
\author[2\authfn{1}]{Axel Cort\'es-Cubero}
\author[1\authfn{2}]{Israel Almod\'ovar-Rivera}
\author[1\authfn{2}]{Wolfgang Rolke}

\contrib[\authfn{1}]{Equally contributing authors.}

\affil[1]{Department of Mathematical Sciences, University of Puerto Rico, Mayag\"uez, Puerto Rico, 00681, USA}
\affil[2]{Protocol Labs, 548 Market Street, San Francisco, CA, 94103, USA}

\corremail{roberto.rivera30@upr.edu}

\fundinginfo{NSF, OAC, Grant/Award Number: 1940179 }

\runningauthor{Rivera R. et al.}

\newtheorem{theorem}{Theorem}
\newtheorem{lemma}{Lemma}

\newtheorem{corollary}{Corollary}

\begin{document}

\maketitle

\section*{Abstract}
In practice, the use of rounding is ubiquitous. Although researchers have looked at the implications of rounding continuous random variables, rounding may also be applied to functions of discrete random variables. For example, to infer the number of excess deaths due to falls after a national emergency, authorities may only provide a rounded average of deaths before and after the emergency started. Deaths from falling tend to be relatively low in most places, and such rounding may seriously affect inference on the change in the rate of deaths. In this paper, we study drawing inference on a parameter from the probability mass function of a non-negative discrete random variable $Y$, when for rounding coarsening width $h$ we get $U=h[Y/h]$ as a proxy for $Y$. We show that the probability generating function of $U$, $\operatorname{E}(U)$, and $\operatorname{Var}(U)$ capture the effect of the coarsening of the support of $Y$. Theoretical properties are explored further under some probability distributions. Moreover, we introduce two relative risks of rounding metrics to aid the numerical assessment of how sensitive the results may be to rounding. Under certain conditions, rounding has little impact. However, we also find scenarios where rounding can significantly affect statistical inference. The methods are applied to infer the probability of success of a binomial distribution and estimate the excess deaths due to Hurricane Maria. The simple methods we propose can partially counter rounding error effects.



\section{Introduction}\label{sec1}
Rounding or rounding off typically replaces a number with an approximated value that is hoped to simplify the representation. Rounding is commonly applied in the life sciences. For example, when measuring occipitofrontal circumference in children (which are important markers of cerebral development), it is recorded to the nearest centimeter, and such a rounding could mask contrasts \citep{wang2013density}. Weight may be rounded to the nearest pound, and age rounded to the nearest year. The effects of rounding on the first two moments of the probability distribution of a continuous random variable have been previously studied \cite{tricker84}, as well as the effects of rounding errors on Type I errors, power, and R
charts \cite{tricker1990effect,tricker1998effect}.  The characteristic function, moments, and oscillatory behavior of rounded continuous random variables were investigated by \cite{janson06}; \cite{paceetal04} studied the properties of likelihood procedures after decimal point rounding, \cite{wang2013density} suggested rounding errors may affect statistical inference, and \cite{chen2021non} defined non-asymptotic moment bounds for rounded random variables.  Many of these studies have concluded that rounded random variables can have properties similar to those of true (hidden) random variable counterparts. However, it is unclear how generally good the approximation is. Moreover,  exponential growth in data \citep{beath2012finding, rivera2019incorporating, rivera2020principles}, recent tendencies in deep learning to lower precision \citep{rodriguez2018lower,wang2018training,colangelo2018exploration,gupta2015deep}, and development of physically informed machine learning models \citep{raissi2017physics,rao2020physics,hooten2011assessing,wikle2010general} make it imperative to better understand the effects of rounding and truncation error \citep{kutz2013data}. 

\begin{figure}[h]
	\begin{center}
		\includegraphics[width=\textwidth]{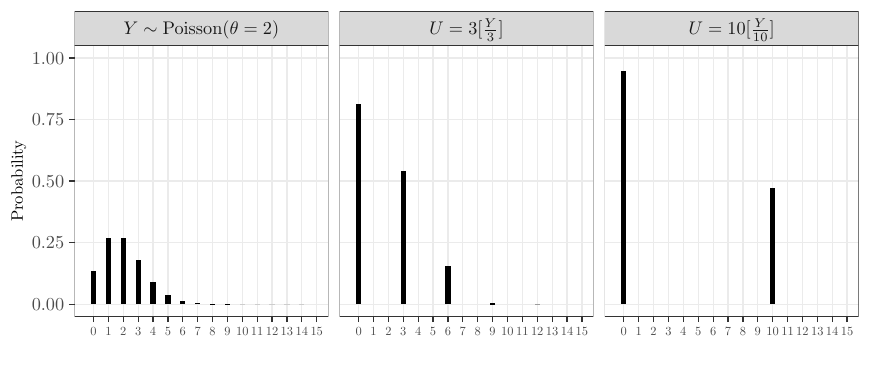}
	\end{center}
	\caption{Left panel shows the pmf of $Y \sim \mbox{Poisson}(\theta=2)$, middle panel presents the pmf of $U$ when $n=3$ while right panel shows the pmf of $U$ when $n=10$.} \label{fig:pmf}
\end{figure}

Discrete data in public health, education, and demography are commonly grouped into intervals, transforming it into categorical or ordinal formats to facilitate tabulation, multinomial modeling, and confidentiality protection. Instead of grouping, in this paper we address rounding that produces a coarser resolution for a discrete random variable. Our emphasis is on the effects of rounding for a non-negative discrete count random variable. Let $Y = X_1 + ... + X_n$ be the total counts from $n$ independent measurements, with some probability mass function (pmf), $P_Y(y)$, parameterized by $\nu$. However, instead of obtaining $Y$ directly, only an average over some measurements rounded to the nearest non-negative integer is available, which is then used to estimate the counts. Define $[X]$ as $X$ rounded to the nearest integer. If $h$ is the rounding coarsening width - the distance between two adjacent points of the rounding lattice \cite{ahmad2024note}, the rounded average random variable is $[\frac{Y}{h}]$ and a rounding-based perturbation of count $Y$ may be expressed as 
\begin{equation*}
 U=per(Y)=h[\frac{Y}{h}]; h \ge 1. 
\end{equation*}
 Whether $u$ has an upper bound or not depends on the pmf of $Y$. $[\frac{Y}{h}]$ has support $0,1,..$ so $u \in \{0,h,2h,...\}$; a multiplier of $h$. In this paper we focus on the situation when $h=n$. When $n=1$, then $U=y$. Since $n$ is fixed, it is possible that $P(U=y)=0$ for some $y$ although $P(Y=y)\neq 0$. For example, if $n=3$, then $P(n[\frac{Y}{n}]= 10)=0$ even if $P(Y=10)$ is quite different from zero. From the support of $U$, it is clear that a noticeable binning of $Y$ values occurs, and the larger the sample size, the more separated the support values of $U$ become. This is a form of coarsening \citep{taraldsen2011analysis}. Figure \ref{fig:pmf} illustrates how $U$ coarsens the support of $Y \sim \mbox{Poisson}(\theta=2)$; the larger the sample size, the more drastic the coarsening. At $n=10$, treating $U$ as a Poisson random variable would lead us to underestimate $\theta$. Because the average is rounded to the nearest integer, attempting to estimate the total by $U$ and treating it as $Y$ will not adequately account for uncertainty. Our aim is to study how the use of $U$ as a proxy for $Y$ affects inference on $\nu$.

		

  Mortality is often underestimated for pandemics, heatwaves, influenza, natural disasters, and other emergencies, times when accurate mortality estimates are crucial for emergency response \citep{rivera2020excess}. The COVID-19 pandemic made it more evident than ever that determining the death toll of serious emergencies is difficult \citep{rosenbaum2021timeliness}. Excess mortality estimates can yield a complementary assessment of mortality. One way statistical models can estimate excess mortality due to an emergency is by comparing the total deaths before the emergency (period 1), with the total deaths post-emergency (period 2). If excess mortality estimates exceed the official death count from the emergency, the official death count may be an underestimate. Excess death models have shown discrepancies with the official death toll from the COVID-19 pandemic \citep{rivera2020excess}, Hurricane Katrina \citep{kutneretal09,stephensetal07}, Hurricane Maria \citep{riveraandrolke18,rivera2019modeling,santosandhoward18,santos2018differential,kishoreetal18}, heatwaves \citep{canoui2006excess,tong2010excess} and other emergencies. 

 On the morning of September 20, 2017, Hurricane Maria slammed into Puerto Rico with sustained winds of 155 mph and higher gusts. The atmospheric event led to devastation across the island. The effects of the wind and rain from the hurricane were felt for more than 24 hours. Cell phone towers had collapsed. For months, tens of thousands of residents were left without utility services. The severe infrastructure damage hindered proper forensic protocol to determine the causes of death.  More than a month after the natural disaster, the official government figure for Maria's death toll, which included direct and indirect deaths, was 55. The Puerto Rican government resisted making detailed data generally available, and in November 2017, the `Centro de Periodismo Investigativo' (CPI) reported that the average daily number of deaths from September 1 to September 19, 2017 was 82 and that it had increased to 117 from September 20 to September 30 \citep{CPI}. While these rounded averages provided by the local government can be used to determine excess deaths \cite{riveraandrolke18}, instead of having access to actual counts, $X_1$ for period 1, and $X_2$ for period 2, we are instead given $[\frac{X_1}{n_1}]$ and  $[\frac{X_2}{n_2}]$ which we then use to get $U_1$ for period 1, $U_2$ for period 2; values that we rely on to infer on the difference in true mean counts of both periods. How would statistical inference be affected when using $U_1$ and $U_2$ as proxies for $X$ and $Y$ respectively when drawing inference on expected difference on mean mortality?

This article is structured as follows. Section \ref{sec:U} presents some of the theoretical properties of the random variable $U$. The special cases where $Y$ follows a Poisson distribution and when it follows a binomial distribution are also studied. We also introduce relative risk measures that can help assess the effect of rounding and, if necessary, implement numerical likelihood maximization. Section \ref{sec:ex} demonstrates the theory and methods developed through two examples: estimating the probability of success of a binomial random variable and estimating excess deaths due to falls after Hurricane Maria. We summarize our findings and their implications in Section \ref{sec:discuss}. Proofs of all the theoretical results are relayed to the Appendix section. 

\section{Properties of the Proxy Random Variable \texorpdfstring{$U$}{U}}\label{sec:U}

Scientists often round data and then misspecify the probability distribution of the proxy random variable. For example, \cite{tilley2019nearshore} rounds raw catch per unit effort fishing data and then models these data as a Poisson random variable. In our context, the proxy random variable $U$ may have a probability distribution that is significantly different than $Y$.
\begin{lemma}
	\label{theorem:pmfu} If $\lfloor x \rfloor$ maps $x$ to the greatest integer less than or equal to $x$, $\lceil x \rceil$ maps $x$ to the least integer greater than or equal to $x$, $Y$ is a non-negative discrete random variable, and $U=n[\frac{Y}{n}]$, then, 
	\begin{eqnarray}
	P(U=u)&=&\sum_{q=0}^{n-1-g(u)}P\Bigl(Y=h(u)+q+g(u)\Bigr) ;\hspace{0.2cm} u\in \{0,n,2n,\ldots\},\label{eq:pmfU}
	\end{eqnarray}
	where, 
	\begin{eqnarray}
	g(u)=\begin{cases}
	\lfloor \frac{n}{2} \rfloor, & u = 0\label{eq:gub}\\
	0, & u \ge 1
	\end{cases}; \:\:\:\mbox{ and } \:\:\:
	h(u)= \lceil u -\frac{n}{2}\rceil = \begin{cases}
	u-\frac{n}{2}+\frac{1}{2}, & n \text{ is odd}\nonumber\\
	u-\frac{n}{2}, & n \text{ is even}\nonumber
	\end{cases}.
	\end{eqnarray}
\end{lemma}
 Note that when $n=1$, then $P(U=u)=P(Y=u)$. Lemma \ref{theorem:pmfu} assumes that the round half-up tie-breaking rule is used. If $n$ is even, the pmf will depend on the type of tie-breaking rule used (a tie is when the fraction of the average is 0.5). If the round half-to-even rule is used, then it can be shown that (see the Appendix), 
\begin{eqnarray}
P(U=u)
&=&\begin{cases}
\sum_{q=0}^{n}P(Y=u-\frac{n}{2}+q), & u/n \text{ is even}\nonumber\\
\sum_{q=0}^{n-2}P(Y=u-\frac{n}{2}+1+q), & u/n \text{ is odd}\nonumber\\
\end{cases}.
\end{eqnarray}
The rest of this paper proceeds according to the round half-up tie-breaking rule. This was a pragmatic choice, as the rule made theoretical results more compact and did not have an effect on the overall conclusions of the paper.  


Turning to moments, observe from (\ref{eq:pmfU}) that the pmf of $U$, for the most part, aggregates probabilities of $Y$. Thus, to express moments of $U$ as a function of moments of $Y,$ a projection is useful. To accomplish this, first, we will derive an expression for the probability generating function (pgf) \citep{resnick1992adventures} of $U$ from the pgf of $Y,$
\begin{equation*}
G_Y(s)=\operatorname{E}(s^Y) = \sum_{y=0}^{\infty}p_ys^y,
\end{equation*}
where $p_y=P(Y=y)$ and the sum converges for any $s \in {\rm I\!R}$ such that $|s|\le  1$.

\begin{theorem}
	\label{theorem:pgfu} Suppose that $Y$ is a non-negative discrete random variable, let $U=n[\frac{Y}{n}]$ and $\omega=\exp(\frac{2\pi i}{n})$. Then, the pgf of $U$ is,
	\begin{eqnarray}
	G_U(s)&=&
	\frac{(s^n-1)}{ns^{n/2-r}}\sum_{j=0}^{n-1}a(j)\frac{G_Y(s/\omega^{j})}{s-\omega^{j}};\label{eq:pgfU} 
	\end{eqnarray}
where 
	\begin{eqnarray}
	r=\begin{cases}
	1, & n \text{ is even}\nonumber\\
	1/2, & n \text{ is odd}\nonumber
	\end{cases};\:\:\:\mbox{ and }\:\:\:
	a(j)=\begin{cases}
	(-1)^{j}, & n \text{ is even}\nonumber\\
(-1)^{j}\omega^{j/2}, & n \text{ is odd}\nonumber
	\end{cases}.
	\end{eqnarray}
\end{theorem}
Some special cases arise when $n=1$, i.e., no rounding effect, when $G_U(s)=G_Y(s)$, and when $n=2$, $G_U(s) =\frac{1}{2}\Bigl((s+1)G_Y(s)-(s-1)G_Y(-s)\Bigr)$.
Notice from (\ref{eq:pgfU}) that $\omega$ combined with $a(j)$ will lead to non-negligible oscillatory behavior of the moments of $U$ as we will see later on. For large $n$, if $p_y \rightarrow 0$ as $y$ increases, then fewer terms in the summation in (\ref{eq:pgfU}) will be different from zero, giving $G_U(s)$ a simpler form. Theorem \ref{theorem:pgfu} helps us find expressions for moments of $U$ as a function of moments of $Y$ and, therefore, better understand the impact of rounding.
\begin{theorem}
\label{theorem:expofU} For any non-negative discrete random variable $Y,$ if $U=n[\frac{Y}{n}]$ and $\omega=\exp(\frac{2\pi i}{n})$, then: 
\begin{equation}
\operatorname{E}(U) = 
\operatorname{E}(Y) + \frac{1}{2} \left(2r-1\right)+\sum^{n-1}_{j=1}a(j)\frac{G_{Y}(\frac{1}{\omega^j})}{1-\omega^j};\label{eq:expofU}
\end{equation}
	and 
\begin{eqnarray}
\operatorname{Var}(U) &=& \operatorname{Var}(Y)+\frac{1}{12}(n^2-1)  -\left(2\operatorname{E}(Y)-1\right)\sum^{n-1}_{j=1}a(j)\frac{G_{Y}(\frac{1}{\omega^j})}{1-\omega^j}\nonumber\\
&& -\left(\sum^{n-1}_{j=1}a(j)\frac{G_{Y}(\frac{1}{\omega^j})}{1-\omega^j}\right)^2+2\sum^{n-1}_{j=1}a(j)\left(\frac{G_{Y}^{'}(\frac{1}{\omega^j})}{\omega^j(1-\omega^j)} -\frac{G_{Y}(\frac{1}{\omega^j})}{\left(1-\omega^j\right)^2}\right).\label{eq:varofU}
\end{eqnarray}
where $r$ and $a(j)$ are defined in Theorem \ref{theorem:pgfu}.
\end{theorem}

Let us examine the results until now and explore the properties of $U$ as $n \to \infty$. According to Lemma \ref{theorem:pmfu}, $P(U=u) \approx P(u -\frac{n}{2} \le Y \le u + \frac{n}{2})$. As seen previously, as $n \to \infty$, the support of $U$ becomes more spread out, and the probability mass of $U$ must concentrate on fewer values of the random variable. 

Equation (\ref{eq:varofU}) is similar to the proposed Sheppard's correction \cite{sheppard1898calculation, tricker1998effect, schneeweissetal10} except that  Sheppard's correction ignores that the rounded random variable $U$ and the rounding error depend on $Y$  \cite{zhao2020bayesian}. Furthermore, (\ref{eq:expofU}) and (\ref{eq:varofU}) include alternating series terms dependent on $\omega$, and for large $n$ the difference between successive terms in each series is small. However, the summation terms in (\ref{eq:expofU}) and (\ref{eq:varofU}) will also depend on distribution parameter values relative to $n$. While Theorem \ref{theorem:expofU} offers insight into the analytic structure of rounding bias, $\omega^j$ and $1-\omega^j$ can be numerically unstable when $n$ is large enough, inhibiting practical application.

\subsection{Poisson Case}
In this section, we will explore the theoretical consequences of the random variable $U$ when $Y \sim \mbox{Poisson}(\theta)$. We aim to estimate $\theta$. From (\ref{eq:pmfU}) the pmf is given by,
\begin{eqnarray}
P(U=u)&=&\sum_{q=0}^{n-1-g(u)}\frac{\theta^{h(u)+q+g(u)}e^{-\theta}}{(h(u)+q+g(u))!} ; \:\: u\in \{0,n,2n,\ldots\}\label{eq:poispmf}\\
&=&\begin{cases}
\theta^{h(u)+g(u)}e^{-\theta}\sum_{q=0}^{n-1-g(u)}\frac{\theta^{q}}{(h(u)+q+g(u))!}; &  n\ge2\nonumber\\
\frac{\theta^{u}e^{-\theta}}{u!}; & n=1
\end{cases}.
\end{eqnarray}
As expected, when $n=1$, i.e., no rounding was done, the pmf of $U$ is the pmf of a Poisson random variable. Our intention is to draw inference on $\theta$ using $U$ instead of $Y$. Specifically, we study whether the rounding leads to substantial differences between $\operatorname{E}(U), \operatorname{Var}(U)$, and $\theta$. If we consider that $Y$ is counting events over $n$ units (e.g., over $n$ days) such that each has independent counts $X_k \sim \mbox{Poisson}(\lambda_k)$, then $Y=\sum_{k=1}^{n}X_k \sim \mbox{Poisson}(\sum_{k=1}^{n}\lambda_k)$ and $\operatorname{E}(U)$ would also increase with $n$ since $\theta=\sum_{k=1}^{n}\lambda_k$. Assuming all $X_k$ are identically distributed, as $n \rightarrow \infty$, $\sum_{k=1}^{n}X_k/n$ converges in probability to $\lambda$ by the law of large numbers. Thus, if $\lambda = 0.4$ and $n$ is large, then  $\sum_{k=1}^{n}X_k/n$ should be close to 0.4, which means that $U$ is close to 0, although $\theta = n\lambda$. In contrast, if $\lambda = 0.6$ and $n$ is large, then $\sum_{k=1}^{n}X_k/n$ should be close to 0.6, which means that $U$ is close to $n$, although $\theta = n\lambda$. That is, we can't generally say that for any $\lambda$, $\operatorname{Bias}(U) \rightarrow 0$ as $n \rightarrow \infty$ (the bias of $U$ as an estimator of $n\lambda$). 
\:
\begin{corollary}\label{theorem:expupoissonN} If $Y \sim \mbox{Poisson}(\theta)$, $U=n[\frac{Y}{n}]$, $n>1$, and $\omega=\exp(\frac{2\pi i}{n})$, then 
	\begin{equation}
	\operatorname{E}(U) = \theta+ \frac{1}{2}(2r-1)+e^{-\theta}\sum^{n-1}_{j=1}a(j)\frac{e^{\frac{\theta }{\omega^j}}}{1-\omega^j}\label{eq:expuPoi}
	\end{equation}
	and
	\begin{equation}
	\operatorname{Var}(U)
	=\theta+\frac{1}{12}(n^2-1)-e^{-2\theta}\left(\sum^{n-1}_{j=1}a(j)\frac{e^{\frac{\theta}{\omega^j}}}{1-\omega^j}\right)^2-\sum^{n-1}_{j=1}a(j)\frac{e^{\frac{\theta}{\omega^j}}}{1-\omega^j}\left(\frac{2e^{-\theta}}{1-\omega^j}-e^{-\theta}\right)\label{eq:varuPoi}
	\end{equation}
	where $r$ and $a(j)$ are defined in Theorem \ref{theorem:pgfu}. 
\end{corollary}
See Appendix \ref{sec:euvarpois} for the proof. When $n=1$, clearly from (\ref{eq:poispmf}), $\operatorname{E}(U) = \operatorname{Var}(U) = \theta$. As stated above, the expressions for the moments of $U$ include alternating series terms that will depend on $\theta$ and $n$. For small $\theta$, $e^{\theta} \approx 1+\theta$. When $n$ is even, for some small values of $\theta$, we see from (\ref{eq:expuPoi}) that $U$ displays substantial bias in estimating $\theta$. The oscillatory behavior of $\operatorname{E}(U)$ and $\operatorname{Var}(U)$ as $\theta$ and $n$ vary is generally not negligible. Specifically, if $n=2$ then 
	\begin{equation*}
	\operatorname{E}(U) = \theta+\frac{1}{2}-\frac{e^{-2\theta}}{2}; \:\:
	\operatorname{Var}(U)
	=\theta+\frac{1}{4}- \frac{e^{-4\theta}}{4}.
	\end{equation*}
 Thus, when $n=2$, $\operatorname{E}(U) \to 0 $ as $\theta \rightarrow 0$. Same for $\operatorname{Var}(U)$. In contrast, if $\theta=0.1$, $\operatorname{E}(U)$ and $\operatorname{Var}(U)$ are approximately twice as large. But what happens to $\operatorname{E}(U)$ and $\operatorname{Var}(U)$ as $\theta$ becomes large? 


\begin{lemma}\label{poissonlarge}
 Let $X_k$ be independent Poisson with parameter $\lambda$, and let $Y=\sum_{k=1}^{n}X_k \sim \mbox{Poisson}(\theta)$, then $\theta = n\lambda$ and for fixed $n$, 
\begin{equation}
    \lim_{\lambda\to\infty}\frac{1}{n\lambda}\operatorname{E}(U)= \lim_{\lambda\to\infty}\frac{1}{n\lambda}\operatorname{Var}(U)=1\nonumber
\end{equation}
\end{lemma}
This result makes intuitive sense. When $\lambda$ is very large relative to $n$, the effect of rounding is small because its fractional part becomes minor. However, when $n$ is much bigger than $\lambda$, the fractional part of $\lambda$ becomes relevant, and rounding will have an impact. 

\subsubsection{Maximum Likelihood Estimator of \texorpdfstring{$\theta$}{theta}}
If the distribution of $U$ is misspecified to be a Poisson with mean $\theta^o$, then $U$ will be a consistent estimator of $\theta^o$, although the random variable may fail to estimate $\theta$ \citep{white1982maximum}. From Corollary \ref{theorem:expupoissonN}, when $Y\sim \mbox{Poisson}(\theta)$, 
\begin{eqnarray}
Bias_{\theta}\ U=\operatorname{E}(U) - \theta = \frac{1}{2}(2r-1)+e^{-\theta}\sum^{n-1}_{j=1}a(j)\frac{e^{\frac{\theta }{\omega^j}}}{1-\omega^j}\nonumber
\end{eqnarray}
In light of the theoretical properties of $U$, we now turn to the estimation of $\theta$ using the likelihood function given the proxy random variable. 

\begin{res}
	\label{theorem:mlealpha} If $Y \sim \mbox{Poisson}(\theta)$, and $U=n[\frac{Y}{n}]$, then the maximum likelihood estimator (MLE) is
	\begin{eqnarray}
	\hat{\theta}
	&=&\prod_{q \in \mathcal{P}} (\lceil u -\frac{n}{2}\rceil+g(u)+q)^{\frac{1}{m}};\label{eq:mlealpha}
\end{eqnarray}
	where $q=0,...,n-1$ and $\mathcal{P}$ is the set such that $\lceil u -\frac{n}{2}\rceil+g(u)+q>0$ and $m$ is the length of $\mathcal{P}$.
\end{res}

If $n=2$, then
 
\begin{eqnarray}
\hat{\theta}&=&\begin{cases}
\Bigl(u(u-1)\Bigr)^{1/2}, & u\ge1\nonumber\\
0, & u=0
\end{cases}
\end{eqnarray}
For even $n$, $\hat{\theta}<u$ except when $u=0$, then $\hat{\theta}=0$ 
. Next, assume independent $X_k \sim \mbox{Poisson}(\lambda)$, and $Y=\sum_{k=1}^{n}X_k \sim \mbox{Poisson}(\theta)$ where $\theta=n\lambda$. Then we obtain the following three theoretical properties of $\hat{\theta}$.

\begin{res}\label{eq:PoissonMLEthetalarge}
For fixed $n$,
\begin{equation}
\lim_{\lambda\to\infty}\frac{1}{n\lambda}\operatorname{E}(\hat\theta)=1\nonumber
\end{equation}
\end{res}

\begin{res}\label{eq:PoissonMLEnlarge}
Let $\lambda>0.5$ and $\lambda\neq I+0.5$, where $I$ is a positive integer. 
If $\hat\theta$ is the MLE of $\theta$, then 
\begin{equation}
\lim_{n\to\infty}\frac{1}{n}\operatorname{E}(\hat\theta)= \left(v_0-\frac{1}{2}\right) e^{\Biggl(\Bigl(v_0+\frac{1}{2}\Bigr)\log\left(\frac{v_0+\frac{1}{2}}{v_0-\frac{1}{2}}\right)-1\Biggr)}\label{largenmse}
\end{equation}
where $v_0=\lfloor \lambda+0.5\rfloor$.

\end{res}
When $\lambda<0.5$, then as $n \to \infty$, $U \to 0$, and $\hat{\theta} > 0$.

\begin{res}\label{eq:PoissonMLEnlargev2}
    
If$\lambda<0.5$, then 
\begin{eqnarray}
\lim_{n\to\infty}\frac{1}{n}\operatorname{E}(\hat\theta)&=& \frac{1}{2e}\nonumber
\end{eqnarray}

\end{res}
Result \ref{eq:PoissonMLEthetalarge} presents conditions for the MLE to be unbiased when $\lambda$ is much larger than $n$, while Result \ref{eq:PoissonMLEnlargev2} states that as $\lambda$ becomes much smaller than $n$, the MLE reaches a bias of at most 0.32
 . Result \ref{eq:PoissonMLEnlarge} indicates that when $n$ is larger than $\lambda >0.5$, the MLE may have a small bias. However, the percentage expected error of the MLE is very small if the fixed lambda is large enough. Also, when $\lambda=I+0.5$, the expected value of the MLE becomes the average of the expected value formulas we just derived for $v_0=I$ and $v_0=I+1$. For $I >0$ we have



\begin{equation*}
\lim_{n\to\infty}\frac{1}{n}\operatorname{E}(\hat\theta)=\frac{1}{2}\left(I-\frac{1}{2}\right) e^{\Biggl(\Bigl(I+\frac{1}{2}\Bigr)\log\left(\frac{I+\frac{1}{2}}{I-\frac{1}{2}}\right)-1\Biggr)}
+\frac{1}{2}\left(I+\frac{1}{2}\right) e^{\Biggl(\Bigl(I+\frac{3}{2}\Bigr)\log\left(\frac{I+\frac{3}{2}}{I+\frac{1}{2}}\right)-1\Biggr)}
\end{equation*}


\subsection{Binomial Case}\label{sec:bm}
We now consider the case where there are independent random variables $X_k\sim \mbox{Binomial}(m,\phi)$, $Y=\sum_{k=1}^{n}X_k$, and the goal is to infer on $\phi$. Now $U$ has domain $u \in \{0,n,...,mn\}$. Once more, $U$ effectively binned the possible values of $Y$. 
\begin{corollary}\label{theorem:binomU} If $Y \sim \mbox{Binomial}(mn,\phi)$, $q=1-\phi$, $U=n[\frac{Y}{n}], n>1$ and $\omega=\exp(\frac{2\pi i}{n})$, then 
	\begin{eqnarray}
\operatorname{E}(U) &=& mn\phi+\frac{1}{2}\left(2r-1\right)+\sum^{n-1}_{j=1}a(j)\frac{(1-\phi+\frac{\phi}{\omega^j})^{mn}}{1-\omega^j};\nonumber
	\end{eqnarray}
	and
	\begin{eqnarray}
	\operatorname{Var}(U)
	&=&mn\phi q+ \frac{1}{12}(n^2-1)
 -\left(2mn\phi-1\right) \sum^{n-1}_{j=1}a(j)\frac{(q+\frac{\phi}{\omega^j})^{mn}}{1-\omega^j}\nonumber\\
 &&-\left(\sum^{n-1}_{j=1}a(j)\frac{(q+\frac{\phi}{\omega^j})^{mn}}{1-\omega^j}\right)^2- 2 \sum^{n-1}_{j=1}a(j)\left(\frac{mn\phi(q+\frac{\phi}{\omega^j})^{mn-1}}{1-\omega^j} -\frac{(q+\frac{\phi}{\omega^j})^{mn}}{\left(1-\omega^j\right)^2}\right).\nonumber
\end{eqnarray}
where $r$ and $a(j)$ are defined in Theorem \ref{theorem:pgfu}. 
\end{corollary}
If $n=2$ then,
\begin{eqnarray}
\operatorname{E}(U) = 2m{\phi}+\dfrac{1}{2}-\frac{\left(1-2{\phi}\right)^{2m}}{2}.\nonumber
\end{eqnarray}
Furthermore, when $m=1$ in this scenario, $\operatorname{E}(U)=2\phi(2- \phi)$ instead of $2\phi = \operatorname{E}(Y)$. Therefore, using $U$ as a proxy to perform inference on $\phi$ would result in substantial bias
. The MLE of $\phi$ in terms of $U$ appears to have a complicated mathematical form. However, it can be determined numerically, as described in the next section.

\subsection{Relative Risk Metrics of Rounding}\label{sec:nummle}

A simple way to compare the risk of two statistics is through the ratio of their mean squared error (MSE).
For parameter $\nu$ we have, 
\[\psi(\nu, n)=\frac{\operatorname{MSE}\left(T_1(Y)\right)}{\operatorname{MSE}\left(T_2(Y)\right)}=\frac{\operatorname{E}\left(T_1(Y)-\nu\right)^2}{\operatorname{E}\left(T_2(Y)-\nu\right)^2}.\]
If the MLE using $y$ is $\hat{\nu}(y)$ and the MLE using $u$ is $\hat{\nu}(u(y,n))$, we propose two relative risks of rounding (RRoR) metrics:

\begin{eqnarray}
\psi_M(\nu, n)=\frac{\operatorname{MSE}\left(\hat{\nu}(u(Y,n))\right)}{\operatorname{MSE}\left(\hat{\nu}(Y)\right)}
	;\:\:\:\mbox{ and }\:\:\:
\psi_U(\nu, n)=\frac{\operatorname{MSE}\left(u(Y,n)\right)}{\operatorname{MSE}\left(\hat{\nu}(u(Y,n))\right)}\nonumber
\end{eqnarray}
Because of the discrete nature of the distributions under consideration, we provide a framework to obtain $\hat{\nu}(u(Y,n))$ numerically when the underlying true count follows a binomial or negative binomial distribution. We proceed as follows (see R routine \emph{RoRR} in supplementary material): 

\begin{algorithm}
\scriptsize
\setstretch{1.35}
\caption{\textsc{Calculate RRoR Metrics}}\label{algorithm:mseratio}
\begin{algorithmic}[1]

\State Select an appropriate set of parameter values for the pmf of $Y$. 
\State Find integers $y$ for which
\(P(Y=y) \le 0.999\) for all parameter values;
\State For values $y$ from the previous step, calculate $u$
\If{$Y \sim \mbox{Poisson}(\theta)$}
\State find the MLE of the parameter via (\ref{eq:mlealpha})
\Else 
\State find the MLE via numerical
optimization of the log-likelihood function using $u$ ;
\EndIf 
\State Using each $y$ and $P(Y=y)$ from Step 1, obtain $\psi_M(\nu,n)$ or $\psi_U(\nu,n)$ by evaluating
\begin{eqnarray}
    \operatorname{E}\left(T(Y)-\nu\right)^2  = \sum_{y=0}^\infty \left(T(y)-\nu\right)^2P(Y=y)\nonumber
\end{eqnarray}
where $T(Y)$ is the appropriate statistic.    
\end{algorithmic}
\end{algorithm}

The output of the routine is either a graph of the selected RRoR or a list of parameter values, mean squared errors and relative risk measures.
\begin{figure}[H]
	\begin{center}
		\includegraphics[width=5in]{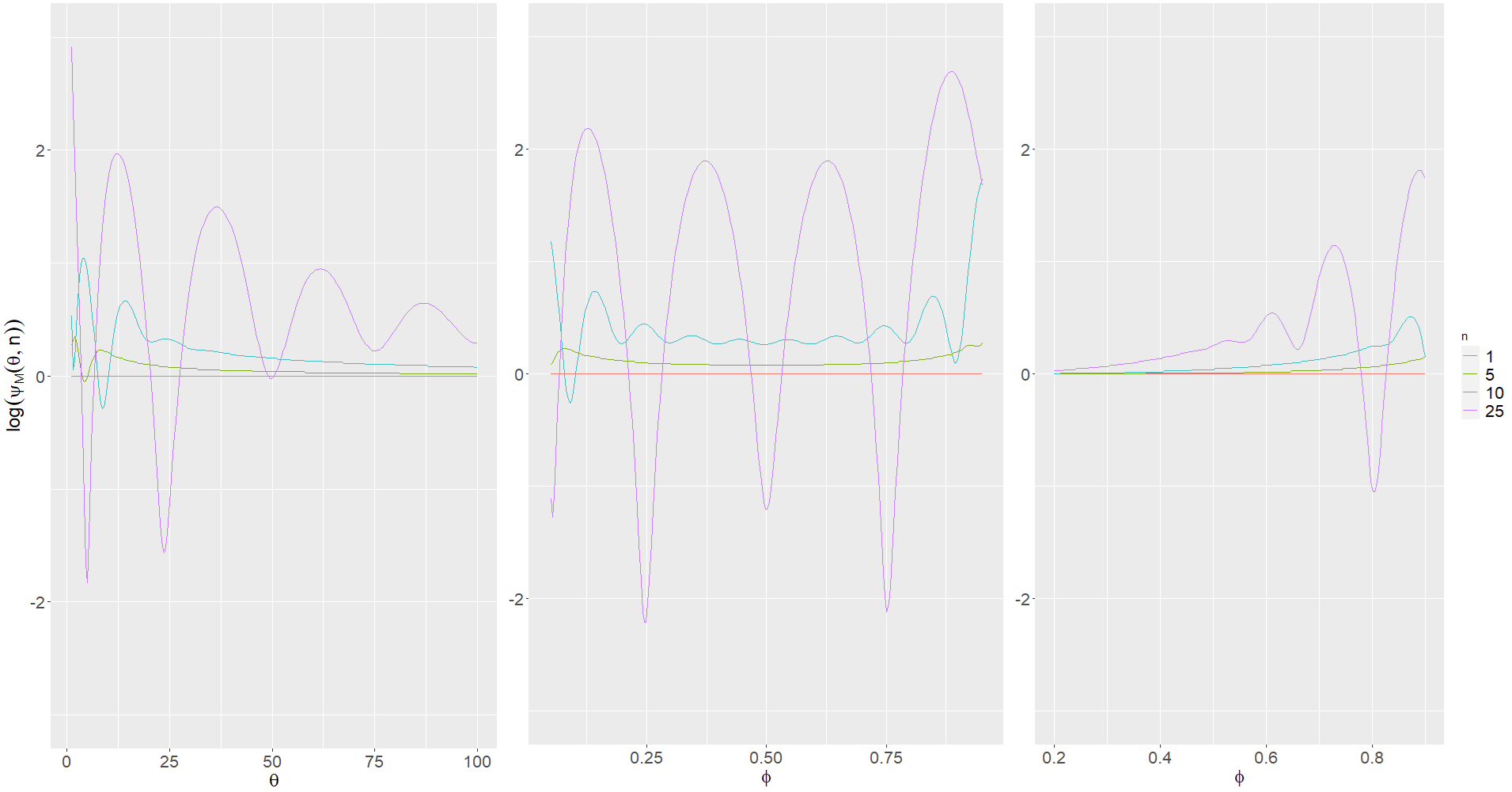}
	\end{center}
	\caption{$\log(\psi_M(\nu, n))$ based on 4 sample sizes when $Y$ follows a Poisson distribution (left), a binomial (middle), and a negative binomial (right).} \label{fig:logmseratio}
\end{figure}

As examples, we consider RRoR metrics in the case of Poisson, binomial, and negative binomial $Y$ when $n = 1, 5, 10, 25$. For the Poisson scenario, our aim is to estimate the expected value $\theta$, while for the other two distributions, the probability of success $\phi$. To visualize results for all sample sizes simultaneously, they are presented on a logarithmic scale. When $Y$ follows a Poisson distribution, $\log(\psi_M(\theta, n))$ is generally larger for small values of $\theta$ and $n$. This makes sense as counts $u$ are biased to small values (left panel, Figure \ref{fig:logmseratio}). Inflection points of $\log(\psi_M(\theta, n))$ occur approximately at $\theta = n/2, n, 3n/2, \ldots$. When $\theta$ is close to $n$, then $\operatorname{MSE}\left(\hat{\theta}(u(Y,n))\right) < \operatorname{MSE}\left(\hat{\theta}(Y)\right)$ because the bias and variance of $\hat{\theta}(u(Y,n))$ become small.  Larger values of $n$ make the range of ratio values larger (and $\operatorname{MSE}(U)$ and $\operatorname{MSE}(\hat{\theta})$ also have a larger range of values) and $log(\psi_M(\theta, n))$ will oscillate at a diminishing amplitude as $\theta$ increases. In the case of $Y$ following a binomial distribution, $\psi_M(\phi, n)$ also oscillates with greater amplitude as $n$ increases (middle panel, Figure \ref{fig:logmseratio}). When $n=25$,  $\operatorname{MSE}\left(\hat{\phi}(u(Y,n))\right) > \operatorname{MSE}\left(\hat{\phi}(Y)\right)$ for some $\phi$ but for other $\phi$ values $\operatorname{MSE}\left(\hat{\phi}(u(Y,n))\right) < \operatorname{MSE}\left(\hat{\phi}(Y)\right)$. Lastly, when $Y$ follows a negative binomial distribution, the right panel of Figure \ref{fig:logmseratio} demonstrates large oscillations in $\log(\psi_M(\phi, n))$, particularly for large values of $\phi$. It should be noted that the results for the binomial and negative binomial will depend on the chosen parameter $N=nm$, and similarly, the Poisson chart will depend on the ratio between $n$ and $\theta$. Nevertheless, this example presents further evidence of rounding effects not being ignorable. Specifically, in some cases, rounding has an effect on estimating the parameter of interest regardless of sample size or parameter value.

\begin{figure}[H]
	\begin{center}
		\includegraphics[width=5in]{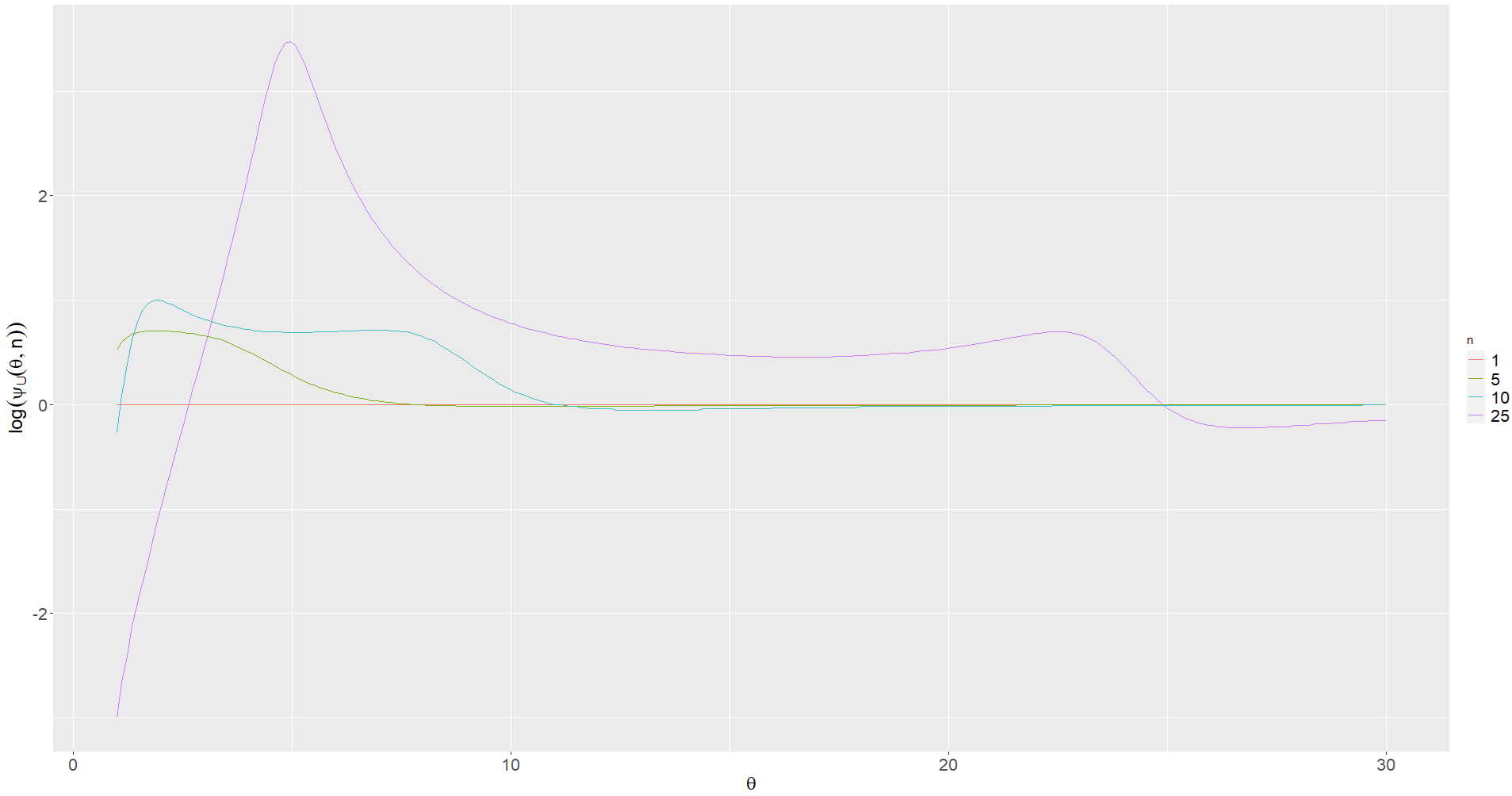}
	\end{center}
	\caption{$\log(\psi_U(\theta, n))$ based on 4 sample sizes when $Y$ follows a Poisson distribution.} \label{fig:logmseratioUforPoisson}
\end{figure}

Results \ref{eq:PoissonMLEthetalarge}, \ref{eq:PoissonMLEnlarge}, and \ref{eq:PoissonMLEnlargev2} showed us some aspects of $\operatorname{E}(\hat{\theta})$, but how does $\hat{\theta}$ generally perform against $U$ as an estimator of $\theta$? To answer this question, we compute $log(\psi_U(\phi, n))$ in the Poisson scenario. We see that $\operatorname{MSE}(\hat{\theta})$ is generally significantly smaller than $\operatorname{MSE}(U)$ until approximately $\theta \ge n$, when  $\operatorname{MSE}(\hat{\theta})$ starts to become slightly larger than $\operatorname{MSE}(U)$ most of the time (Figure \ref{fig:logmseratioUforPoisson}). The exception is when $\theta$ is very small, the MLE struggles to be close to $\theta$ (Result \ref{eq:PoissonMLEnlargev2}). There is evidence of oscillations in mean squared errors, with a dip when $\lambda$ tends to an integer.  
The range of $\log(\psi_U(\theta, n))$ increases as $n$ increases. For $\theta$ values close to 0 and large $n$, $\hat{\theta}$ has a significantly larger mean squared error than $U$. This is because as $n \rightarrow \infty$, $\lceil u -\frac{n}{2}\rceil$ becomes smaller. This can cause the MLE to use as few as $n/2$ terms, which will bias its result such that $\hat{\theta}>\theta$. In practice, by choosing a suitable range of parameter values for the pmf of $Y$, RRoR metrics can enable a better understanding of the impact of rounding on parameter estimation. 

\section{Applications}\label{sec:ex}
In this section, we first consider drawing inference on the probability of success $\phi$ when the hidden random variable $Y$ follows a binomial distribution. Secondly, we present a real data application, when counts are estimated according to two averages rounded to the nearest non-negative integer coming from two separate time periods, and we wish to draw an inference on the difference of the mean total counts. 

\subsection{Inference on probability of success \texorpdfstring{$\phi$}{phi}}
\begin{figure}[H]
	\begin{center}
		\includegraphics[width=5in]{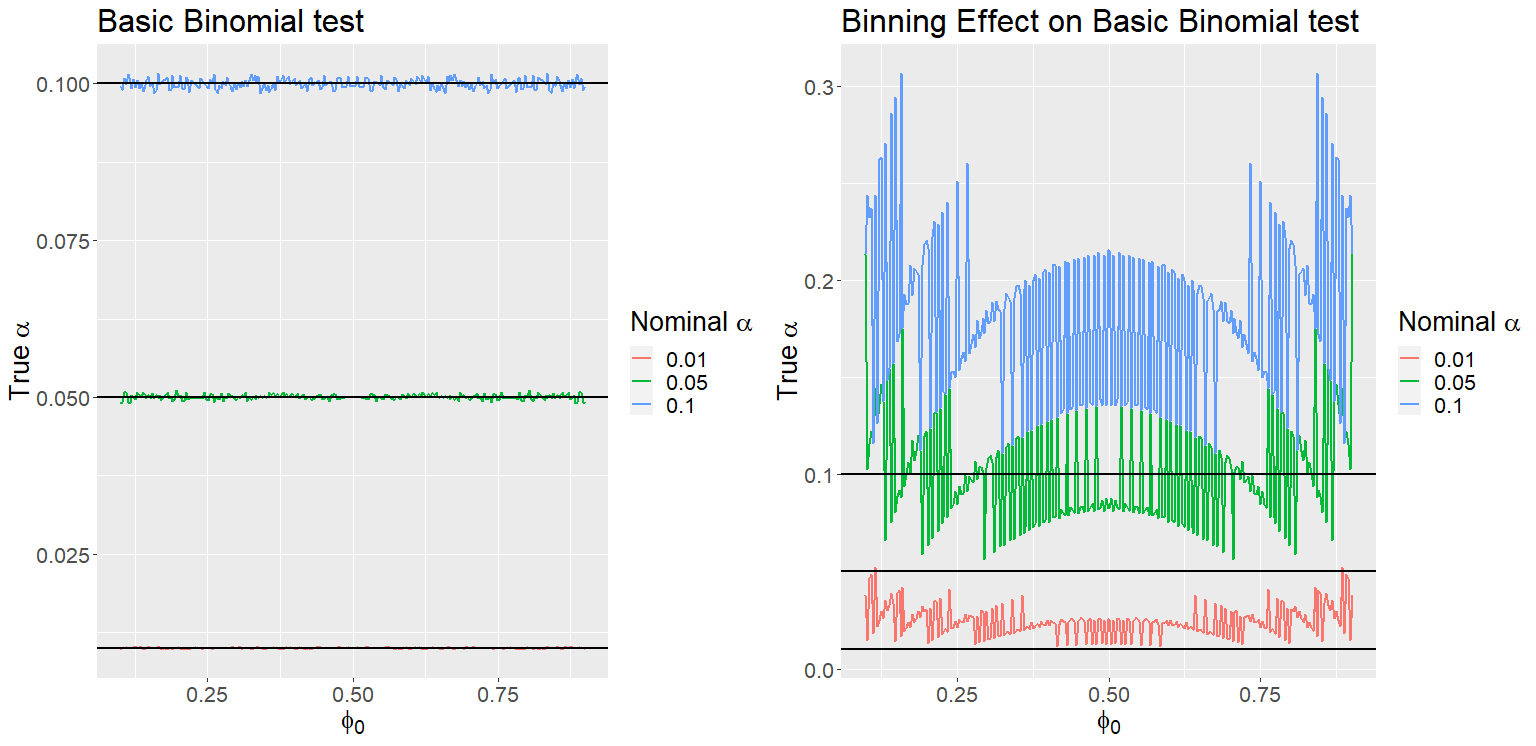}
	\end{center}
	\caption{On the left, we see the true significance level if $Y$ were available. On the right, we have the true significance level when using $U$ and assuming it follows a binomial distribution.} \label{fig:truesig}
\end{figure}

Suppose there is a sequence of latent random variables $X_{k} \sim \mbox{Binomial}(m, \phi)$, $k = 1,\ldots,n$, and our aim is to draw inference on $\phi$. Clearly, $Y = \sum_{k} X_{k} \sim \mbox{Binomial}(mn, \phi)$. Corollary \ref{theorem:binomU} demonstrates how moments of $U$ theoretically deviate from moments of $Y$ and section \ref{sec:bm} gives an example where $\operatorname{E}(U)$, and $\operatorname{E}(Y)$ can be very different. We now examine the practical implications of the theoretical results presented by comparing the true significance level when using $Y$ vs. actually having $U =n[Y/n]$ available. When using $U$, many analysts will draw inferences on $\phi$ by misspecifying its distribution as $\mbox{Binomial}(mn, \phi)$. 
 Specifically, we will test
$$H_o:\phi=\phi_o\text{ vs. }H_a:\phi\ne \phi_o$$ with test statistic $$W=\frac{w-mn\phi_o}{\sqrt{mn\phi_o(1-\phi_o)}}$$
where the null is rejected at significance level $\alpha$ if $\mid W \mid$ is greater than the standardized score $z_{\alpha/2}$ and $w$ is either $y$ or $u$. The true significance level is \citep{casellaandberger01}  $$P(W\le mn\phi_o-z_{\alpha/2}\sqrt{mn(1-\phi_o)}) + P(W\ge mn\phi_o+z_{\alpha/2}\sqrt{mn(1-\phi_o)})$$ To ensure the normal approximation is good we choose $m=500, n=31$ and values of $\phi_o$ between 0.1 and 0.9. Comparison was based on 0.01, 0.05, and 0.1 nominal significance levels. The left panel of Figure \ref{fig:truesig} shows the true significance when $Y$ is available. The oscillatory behavior in true significance can be attributed to the lattice structure in $Y$ \citep{brownetal01}. When using $U$ and misspecifying its distribution as binomial, the true significance levels oscillate as a function of $\phi_o$ much more than when $Y$ is available, with values that can be far higher than the nominal significance level (right panel Figure \ref{fig:truesig}). With $U$, the true significance value is always higher than the nominal $\alpha$. 

\begin{figure}[H]
	\begin{center}
		\includegraphics[width=5in]{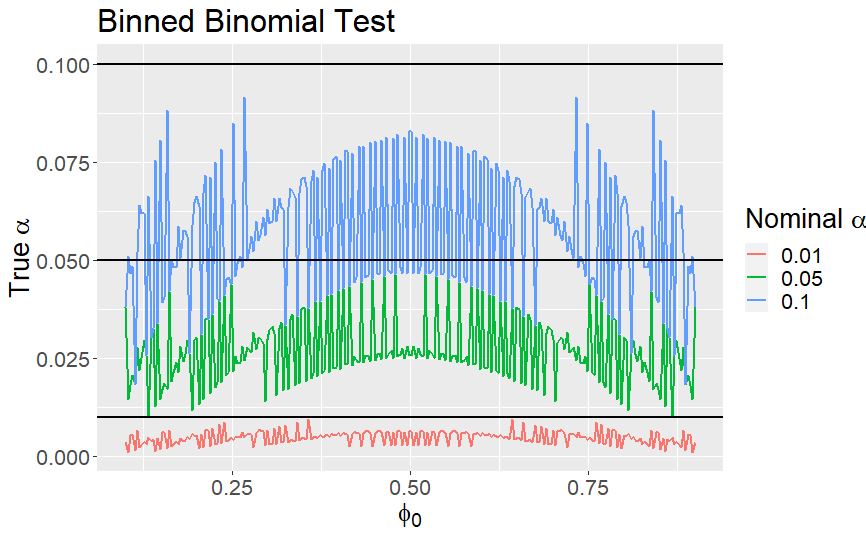}
	\end{center}
	\caption{The true significance level when using $U$ with probability mass function as given in Lemma \ref{theorem:pmfu} or $P(U=u)=\sum_{j=0}^{mn}P(Y=j)I(n[j/n]=u)$ where $Y \sim \mbox{Binomial}(mn, \phi_o)$.} \label{fig:truesigbinnedbinomialtest}
\end{figure}

Instead of misspecifying the distribution of $U$, Figure \ref{fig:truesigbinnedbinomialtest} shows the true significance level when using $U$ where calculations are based on the pmf of $U$ according to Lemma \ref{theorem:pmfu} and where what we call a binned binomial test is performed.
Specifically, note that $P(U=u)=\sum_{j=0}^{mn}P(Y=j)I(n[j/n]=u)$ where $I()$ is an indicator function.  Although the true significance value is now always lower than the nominal $\alpha$, the bias of the true significance level can be much smaller than when misspecifying the distribution of $U$. However, the use of $U$ has caused the oscillations in interval coverage to be much more pronounced in comparison to using $Y$. For example, when the nominal significance level is 0.1, the true significance level of $U$ may be closer to 0.025 for some values $\phi_o$, and when the nominal significance level is 0.05, the true significance level of $U$ may be closer to 0.01 for some values $\phi_o$ (Figure \ref{fig:truesigbinnedbinomialtest}). \textsf{R} code for the binned binomial test is available as supplementary material \citep{rcran20}.

\subsection{Estimating Excess Deaths Due to an Emergency}

Our second application is based on the before-and-after comparison to estimate excess deaths due to an emergency. Daily deaths are assumed to be independent and modeled according to a Poisson distribution with expected value $\lambda$. Then $X \sim \mbox{Poisson}(n_1\lambda)$ represents the total deaths occurring in $n_1$ days before the emergency, and $Y \sim \mbox{Poisson} (n_2\lambda +\beta)$ represents the deaths in $n_2$ days after the emergency starts. $X$ and $Y$ are independent. $\beta$ measures excess deaths, a proxy of the impact of the emergency on mortality. If $X$ and $Y$ were available, a reasonable point estimator of excess deaths would be \citep{riveraandrolke18,rivera2019modeling}

\begin{eqnarray}
   L^* \equiv (\bar{Y} - \bar{X})n_2 = (Y - \frac{n_2}{n_1}X)\label{eq:excdor}
\end{eqnarray}
where $\bar{X} = \frac{X}{n_1}$ and $\bar{Y} = \frac{Y}{n_2}$. For this estimator,
\begin{equation}
  \operatorname{E}(L^*) = \beta \label{eq:norounding}  
\end{equation}
and 
\begin{equation}
 \operatorname{Var}(L^*) = \beta + n_2\left(1+\frac{n_2}{n_1}\right)\label{eq:norounding2}  
\end{equation}
The second term makes adjustments to the variance dependent on the before-and-after emergency sample sizes. However, when $X$ and $Y$ are not available, and total counts must be estimated through averages rounded to the nearest non-negative integer, then the estimator becomes
\begin{eqnarray}
   ([\bar{Y}] - [\bar{X}])n_2 = U_2 - [\bar{X}]n_2 = \left(U_2-\frac{n_2}{n_1}U_1\right) \equiv L\label{eq:excd}
\end{eqnarray}
That is, $U_2-U_1$ where $U_1=n_1[\frac{X}{n_1}]$, is only a suitable estimator when $n_1=n_2$. Our theoretical results shed light on the impact of supplanting (\ref{eq:excdor}) with (\ref{eq:excd}). Referring to Corollary \ref{theorem:expupoissonN} we have, 

\begin{eqnarray}
 \operatorname{E}(L) &=& \beta +\frac{1}{2}(2r_2-1) - \frac{n_2}{2n_1}(2r_1-1) +\nonumber\\ &&e^{-(n_2\lambda+\beta)}\sum^{n_2-1}_{j=1}a_2(j)\frac{e^{\frac{(n_2\lambda+\beta)}{\omega_2^j}}}{1-\omega_2^j}-\frac{n_2 e^{-n_1 \lambda}}{n_1}\sum^{n_1-1}_{j=1}a_1(j)\frac{e^{\frac{n_1 \lambda }{\omega_1^j}}}{1-\omega_1^j}\label{eq:emari}
\end{eqnarray}
where $\omega_k=\exp(\frac{2\pi k}{n_k})$ and

\begin{eqnarray}
	r_k=\begin{cases}
	1, & n_k \text{ is even}\nonumber\\
	1/2, & n_k \text{ is odd}\nonumber
	\end{cases};\:\:\:\mbox{ and }\:\:\:
	a_k(j)=\begin{cases}
	(-1)^{j}, & n_k \text{ is even}\nonumber\\
(-1)^{j}\omega^{j/2}_k, & n_k \text{ is odd}\nonumber
	\end{cases}.
	\end{eqnarray}
Moreover,
\begin{equation}
 \operatorname{Var}(L) = \beta + n_2\left(1+\frac{n_2}{n_1}\right)+\frac{n_2^2 -1}{12} + \frac{n_2^2(n_1^2 -1)}{12n_1^2}+Q_1\label{eq:mari}
\end{equation}
where $Q_1$ is a term resulting from the series in (\ref{eq:varuPoi}). Alternatively, through the invariance property of MLEs \citep{casellaandberger01}, $\widehat{\theta + \beta} - \frac{n_2}{n_1}\hat{\theta}$ can be used as an MLE estimator for $\beta$; where the first term is a function of $U_2$ and the second of $U_1$. Considering the application, it is reasonable to assume that $n_1, n_2$ are not large. The effects of rounding on the expected value of the estimator can be studied comparing (\ref{eq:norounding}) and (\ref{eq:emari}), while rounding effects on estimator variance can be studied comparing (\ref{eq:norounding2}) and (\ref{eq:mari}). We remark:
\begin{itemize}
    \item If $\theta$ is large, then from (\ref{eq:excd}) and Lemma \ref{poissonlarge} we have $\operatorname{E}(L) \approx \operatorname{E}(L^*)$ and $\operatorname{Var}(L) \approx \operatorname{Var}(L^*)$. When $\theta$ or $\theta + \beta$ are large, their respective MLEs $\hat{\theta}$, $\widehat{\theta+\beta}$ should perform well (Lemma \ref{poissonlarge}).
    \item If $\theta$ is not large, $n_1$ is even and $n_2 > 2n_1$, then from (\ref{eq:emari}) we see that $\operatorname{E}(L)$ will deviate considerably from $\operatorname{E}(L^*)$. When $n_2 > n_1$, (\ref{eq:mari}) shows that $\operatorname{Var}(L)$ will deviate considerably from $\operatorname{Var}(L^*)$ regardless of whether $n_1$ is even or odd. Moderate values of $n_1$ and $n_2$ would create a bias due to the third and fourth term in (\ref{eq:mari}). The level of the bias is dependent on $\beta$ and $\theta$, which impact $Q_1$.
    
    \item As implied in section \ref{sec:U}, if either $\theta$ or $\theta + \beta$ are of form $n_k(I+0.5)$, then $\operatorname{Var}(L)$ will be large. Both parameters having this form will result in a larger value of (\ref{eq:mari}).
    
    \item  Corollary \ref{theorem:expupoissonN} and Figure \ref{fig:logmseratioUforPoisson} imply that for $n_2/4 \le \theta + \beta \le n_2/2$, $U_2$ will have a substantially larger $\operatorname{MSE}$ than the MLE of $\theta + \beta$. If $n_1 = n_2$ and $n_1/4 \le \theta \le n_1/2$, $U_1$ will have a substantially larger $\operatorname{MSE}$ than the MLE of $\theta$. This would lead to an overestimation of excess deaths unless $n_2/4 \le \theta + \beta \le n_2/2$; when an underestimation may occur.
\end{itemize}
\subsubsection{Estimating Excess Deaths Due to Falls After Hurricane Maria}
The analysis above demonstrates the theoretical effects rounding would have on estimating excess deaths. We now present a more concrete application. After landfall of Hurricane Maria, the Puerto Rico Vital Statistics Office refused to share death certificate data with researchers or journalists for months,
and in June 2018, a judge ordered the agency to do so. Here, we consider daily death data due to a fall and assume that they follow a Poisson distribution. Given the data, we  illustrate how an analysis based on rounded average daily deaths can lead to problems. Assessing excess deaths by cause can be very valuable in the development of emergency plans \citep{lugo2023closer}. To estimate excess deaths due to a fall from September 20, 2017 to October 31, 2017, we make a comparison to a pre-Hurricane Maria period of the same number of days: August 8  2017 to September 19, 2017. There were 24 deaths caused by a fall in the pre-Maria period, 36 in the post-Maria period which corresponds to 0.57 and 0.86 average daily deaths respectively.  From (\ref{eq:excd}) we have 0 estimated excess deaths, although from (\ref{eq:excdor}) 12 excess deaths are obtained. Clearly, the MLE counterpart, $\widehat{\theta + \beta}-\frac{n_2}{n_1}\widehat{\theta}$, is 0 since $\widehat{\theta + \beta}$ and $\widehat{\theta}$ are both based on $U=42$. 


\begin{figure}[H]
	\begin{center}
		\includegraphics[width=5in]{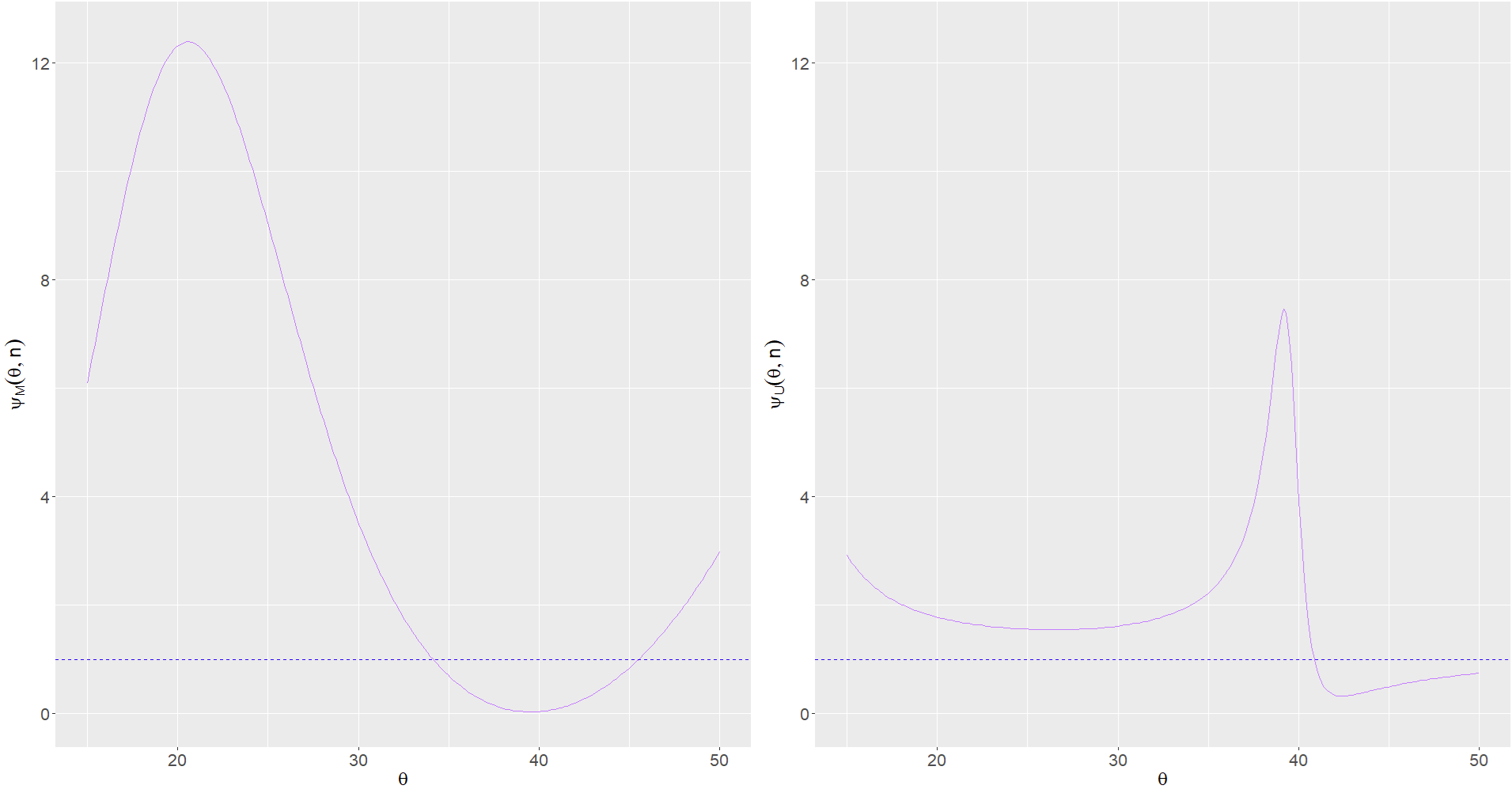}
	\end{center}
	\caption{$\psi_M(\theta, n)$ (left) and $\psi_U(\theta, n)$ (right) when $Y$ follows a Poisson distribution and $n=42$. Horizonal lines represent $\psi_M(\theta, n)=\psi_U(\theta, n)=1$.} \label{fig:mseratiosbyexcdeathsbyfall}
\end{figure}

According to the summary statistics of the available data, we can reasonably expect the mean deaths due to a fall over the 42 days before and after Hurricane Maria to be between 15 and 50. Figure \ref{fig:mseratiosbyexcdeathsbyfall} helps us analyze the sensitivity to rounding in our excess deaths due to a fall application. For true parameter values no greater than approximately 34, $\psi_M(\theta, n)$ is largest. If the estimate of 24 deaths due to a fall before Hurricane Maria is close to the true expected value, rounding would cause $\operatorname{MSE}\Big(\hat{\theta}(u(X,42))\Big)$ to be more than 10 times larger than $\operatorname{MSE}\Big(\hat{\theta}(X)\Big)$. In contrast, if the
estimate of 36 deaths due to a fall after Hurricane Maria is relatively close to the true expected value, 
rounding would cause $\operatorname{MSE}\Big(\hat{\theta}(u(Y,42))\Big)$ to be a little less than twice $\operatorname{MSE}\Big(\hat{\theta}(Y)\Big)$. On the other hand, regardless whether we use $X$ or $Y$,     $\operatorname{MSE}\Big(\hat{\theta}(u(\cdot,42))\Big)$ is lower than $\operatorname{MSE}\Big(u(\cdot,42)\Big)$ with the largest discrepancy when $\theta$ is close to 39. For larger parameter values $\psi_U(\theta, n)$ drops quickly below 1, remaining just below that value as $\theta$ continues to increase. Consequently, if we used estimators based on rounded averages, there is a substantial risk of misrepresenting the true excess deaths due to falls.

\section{Discussion}\label{sec:discuss}
The explosion of data and the proposal of lower-precision deep learning algorithms to speed up computations has made scientists rethink ignoring rounding error. In this paper, we study the effects of relying on an average rounded to the nearest non-negative integer times $n$ measurements to obtain $U$ as a proxy of total counts. We derive expressions for $P(U=u), G_U(s), \operatorname{E}(U)$, and $\operatorname{Var}(U)$. As far as we know, this is the first time the effect of rounding is assessed for discrete random variables. Conditions when rounding error is negligible and when it is not, are presented. For a long time, it was considered that rounding had negligible consequences in statistical inference. Yet the alternating series found in $G_U(s)$ and the moments of $U$ can result in an oscillating behavior dependent on $n$ and parameter values, which can significantly alter statistical inference as the two examples demonstrate. Studies have reached similar conclusions assessing the impact of rounding error on continuous variables and statistical inference \citep{wang2013density,tricker1998effect}. As illustrated through the excess deaths example, rounding may result in significant first-order bias as well. Equation (\ref{eq:excd}) combined with the work from \cite{janson06} may elucidate the influence of rounding when comparing the means of two different continuous random variables. 

We demonstrate how the use of the true pmf of $U$ helps reduce the bias in significance level calculations, albeit the bias may still be substantial. We also present a maximum likelihood estimator for the case of $Y \sim Poisson(\theta)$ and explore its theoretical properties. Furthermore, we introduce two RRoR metrics to aid rounding sensitivity assessment. When estimating excess deaths due to falls after Hurricane Maria the MLE $\hat{\theta}$ performs well for most values of $\theta$ and $\operatorname{MSE}(\hat{\theta})$ is generally lower than $\operatorname{MSE}(U)$ for small parameter values. Code to obtain RRoR metrics and the numerical MLE when the underlying distribution is binomial or negative binomial is provided as supplementary material.

We did not explore methods that calibrate rounding errors. A valuable research path is the further development of relative risk measures analogous to $\psi_M(\theta, n)$ and $\psi_U(\theta, n)$ for continuous random variables; in which case an approximation of the mean squared errors of rounded-based variables will be required. Additional future research includes following a Berkson measurement error model, such that a nonparametric estimator of the distribution of $Y$ could be constructed \citep{wang2013density}.  The optimal transport theory approach \citep{peyre2019computational} is another promising research path. Future work could examine the Wasserstein distance between the distributions of $U$ and $Y$, develop minimum Wasserstein distance estimators \citep{bernton2019parameter} or approximate intractable $U$ distributions \citep{torres2021survey}.

\begin{itemize}
\item The authors gratefully acknowledge comments and suggestions from the associate editor and reviewers, which significantly improved the article. Financial support from the National Science Foundation (NSF) OAC Award 1940179 partially supported Cort\'es and Rivera. 
\item The authors have no conflicts of interest to declare that are relevant to the content of this article.
\item Ethics approval `Not applicable'
\item Consent to participate `Not applicable'
\item All authors consent to the publication of this paper.
\item Data will be made available in the corresponding author's GitHub.
\item Simulation code is available upon request. 
\item Authors' contributions Conceptualization: Rivera; Methodology: Rivera, Cort\'es; Formal analysis and investigation: Rivera, Cortes, Almod\'ovar-Rivera, Rolke; Writing - original draft preparation: Rivera, C\'ortes, Almod\'ovar-Rivera ; Writing - review and editing: Rivera, Almod\'ovar-Rivera; Supervision: Rivera.
\end{itemize}



\begin{appendices}

\section{Proofs}\label{secA1}

\subsection{Proof of Lemma \ref{theorem:pmfu}}\label{app:pmfu}
\begin{proof}
	
It is straightforward to show that,

\begin{eqnarray}
P(U=u)=P(n[\frac{Y}{n}]=u)=P([\frac{Y}{n}]=\frac{u}{n})\nonumber
\end{eqnarray}

The pmf of $U$ depends on $n$. Specifically, when $n$ is odd then,
\begin{eqnarray}
P(U=u)&=&P(u-\frac{n}{2}+\frac{1}{2}\le Y\le u+\frac{n}{2}-\frac{1}{2}) \hspace{1cm} \text{Since } Y \text{ must be an integer}\nonumber\\
&=&\sum_{q=g(u)}^{n-1}P(Y=u-\frac{n}{2}+\frac{1}{2}+q) \hspace{2cm} u\in \{0,n,2n,..\}\nonumber
\end{eqnarray}

Assuming $n$ is even, then things get a bit more complicated, mainly because the pmf will depend on the type of tie-breaking rule used. If the round half-to-even rule is used, then 
\begin{eqnarray}
P(U=u)&=&P(\frac{u}{n}-0.5\le \frac{Y}{n}< \frac{u}{n}+0.5)\nonumber\\
&=&\begin{cases}
\sum_{q=0}^{n}P(Y=u-\frac{n}{2}+q) & u/n \text{ is even}\nonumber\\
\sum_{q=0}^{n-2}P(Y=u-\frac{n}{2}+1+q) & u/n \text{ is odd}\nonumber\\
\end{cases}
\end{eqnarray}
where $u\in \{0,n,2n,..\}$. When $u=0, n\ge2$, then $u-\frac{n}{2} +q < 0$, and $P(Y=u-\frac{n}{2} +q)=0$ until $q \ge \frac{n}{2}$. Alternatively, we may use a round half up tie breaking rule,
\begin{eqnarray}
P(U=u)&=&P([\frac{Y}{n}]=\frac{u}{n})=P(\lfloor\frac{Y}{n}+0.5\rfloor=\frac{u}{n})\nonumber\\
&=&P(\frac{u}{n}\le \frac{Y}{n} + 0.5< \frac{u}{n}+1)= P(\frac{u}{n}-0.5\le \frac{Y}{n}< \frac{u}{n}+0.5)\nonumber\\
&=&
P(u-\frac{n}{2}\le Y< u+\frac{n}{2})=\sum_{q=g(u)}^{n-1}P(Y=u-\frac{n}{2}+q)\nonumber
\end{eqnarray}
Adjusting the summation index to start at zero completes the proof.
\end{proof}

\subsection{Proof of Theorem \ref{theorem:pgfu}}
\begin{proof}
For even $n$, 
\begin{eqnarray}
G_U(s)&=& (p_0 + \cdots + p_{\frac{n}{2}-1})+(p_{\frac{n}{2}} + \cdots + p_{n+\frac{n}{2}-1})s^n + (p_{n+\frac{n}{2}} + \cdots + p_{2n+\frac{n}{2}-1})s^{2n}+\cdots\nonumber\\ && +(p_{(m-1)n+\frac{n}{2}} + \cdots + p_{mn+\frac{n}{2}-1})s^{mn}+\cdots\nonumber\\
&=&\sum_{k=0}^{n/2-1}p_{k}+\sum_{l=0}^{\infty}\sum_{k=0}^{n-1}p_{n/2+ln+k}s^{n/2+ln}\label{eq:ln}
\end{eqnarray}
Recall that the sum of the pgf converges for any $s \in {\rm I\!R}$ such that $\mid s \mid \le 1$. Meanwhile, we may write $G_Y(s)$ as 
\begin{eqnarray}
G_Y(s) = \sum_{k=0}^{n/2-1}p_{k}s^k+\sum_{l=0}^{\infty}\sum_{k=0}^{n-1}p_{n/2+ln+k}s^{n/2+ln+k}.\label{eq:def1}
\end{eqnarray}
Next, we transform (\ref{eq:def1}) the following way,
\begin{eqnarray}
\omega^{jn/2}G_{Y}(s/\omega^j) &=&\sum_{k=0}^{n/2-1}p_ks^k\omega^{j(n/2-k)}+ \sum_{l=0}^{\infty}\sum_{k=0}^{n-1}p_{n/2+ln+k}s^{n/2+ln+k}\omega^{-j(ln+k)} \nonumber\\
&=& \sum_{k=0}^{n/2-1}p_ks^k\omega^{j(n/2-k)}+ \sum_{l=0}^{\infty}\sum_{k=0}^{n-1}p_{n/2+ln+k}s^{n/2+ln+k}\omega^{-jk}\label{eq:transf}
\end{eqnarray}
where the second equality is due to $\omega^{-jln}=1$, for integer values of $j$. The inverse discrete Fourier transform of this function is
\begin{eqnarray}
\frac{1}{n}\sum_{j=0}^{n-1}\omega^{j\left(\frac{n}{2}+q\right)}G_Y\left(\frac{s}{\omega^j}\right)&=&\frac{1}{n}\sum_{j=0}^{n-1}(-1)^j\omega^{jq}G_Y\left(\frac{s}{\omega^j}\right)\nonumber\\
&=&p_{q-n/2}s^{q-n/2}+\sum_{l=0}^\infty p_{\frac{n}{2}+ln+q}s^{\frac{n}{2}+ln+q}\label{eq:idft},
\end{eqnarray}
where $p_k=0$, for any $k<0$.

The probability generating function for $U$ can then be written as,
\begin{eqnarray}
G_U(s)&=&\sum_{q=0}^{n/2-1}p_{q}+\sum_{l=0}^\infty\sum_{q=0}^{n-1}p_{n/2+ln+q}s^{n+ln}\nonumber\\
&=&\sum_{q=0}^{n-1}\left(s^{n/2-q}\frac{1}{n}\sum_{j=0}^{n-1}\omega^{jq}(-1)^jG_Y(s/\omega^j)\right)\nonumber\\
&=& \frac{s^n-1}{n s^{n/2-1}}\sum_{j=0}^{n-1}(-1)^j\frac{G_Y\left(\frac{s}{\omega^j}\right)}{s-\omega^j},\label{pgfuu}
\end{eqnarray}
where the last equality follows from resuming the $q$-dependent terms as a geometric series.

For odd $n$ 
\begin{eqnarray}
G_U(s)=\operatorname{E}(s^U) = (p_0 + \cdots + p_{\frac{n}{2}-\frac{1}{2}})+(p_{\frac{n}{2}+\frac{1}{2}} + \cdots + p_{n+\frac{n}{2}-\frac{1}{2}})s^n + \cdots\nonumber
\end{eqnarray}
and following a similar procedure as for $n$ even we get,
\begin{eqnarray}
G_U(s)&=&\sum_{q=0}^{n/2-1/2}p_{q}+\sum_{l=0}^\infty\sum_{q=0}^{n-1}p_{n/2+1/2+ln+q}s^{n+ln}\label{eq:ln2}\\
&=&\sum_{q=0}^{n-1}\left(s^{n/2+1/2-q}\frac{1}{n}\sum_{j=0}^{n-1}\omega^{j(q+1/2)}(-1)^jG_Y(s/\omega^j)\right)\nonumber\\
&=&\frac{s^n-1}{n s^{(n-1)/2}}\sum_{j=0}^{n-1}(-1)^j\frac{\omega^{j/2}G_Y\left(\frac{s}{\omega^j}\right)}{s-\omega^j}\nonumber
\end{eqnarray}
\end{proof}

\subsection{Proof of Theorem \ref{theorem:expofU}}
Starting with the expected value, we will consider the version of $G_U(s)$ free of pole singularity at $s=1$; thus for $j=0$ in (\ref{eq:pgfU}) we have:
\begin{eqnarray}
\frac{s^n-1}
{s-1}\nonumber 
\end{eqnarray}
a finite geometric series (when $s \neq 1$) and therefore:
\begin{eqnarray}
G_U(s) &=&\frac{1}{n}\left(G_Y(s)\sum_{t=0}^{n-1}s^{t+r-\frac{n}{2}}+\frac{s^n-1}{s^{\frac{n}{2}-r}} \sum^{n-1}_{j=1}a(j)\frac{G_Y(s/\omega^j)}{s-\omega^j}\right)\nonumber
\end{eqnarray}

Taking the derivative with respect to $s$ we have:

\begin{eqnarray}
n G'_U(s) &=& \frac{G_{Y}(s)}{2}\sum_{t=1}^{n}(2r-2+2t-n)s^{r-2+t-\frac{n}{2}}+G_{Y}^{'}(s)\sum_{t=1}^{n}s^{r-1+t-\frac{n}{2}} +\nonumber\\
&& \frac{1}{2}\Bigl((n+2r)s^{\frac{n}{2}+r-1}-(-n+2r)s^{-\frac{n}{2}+r-1}\Bigr) \sum^{n-1}_{j=1}a(j)\frac{G_{Y}(\frac{s}{\omega^j})}{s-\omega^j}+\nonumber\\
&& \left(s^{\frac{n+2r}{2}}-s^{\frac{-n+2r}{2}}\right) \sum^{n-1}_{j=1}a(j)\left(\frac{G_{Y}^{'}(\frac{s}{\omega^j})}{\omega^j(s-\omega^j)} -\frac{G_{Y}(\frac{s}{\omega^j})}{\left(s-\omega^j\right)^2}\right)\label{eq:firstderivpgf}
\end{eqnarray}

Thus, at $s=1$ we get:

\begin{eqnarray}
n G'_U(1) &=& \frac{1}{2}G_{Y}(1) \left(2r-1\right)n+nG_{Y}^{'}(1)+n\sum^{n-1}_{j=1}a(j)\frac{G_{Y}(\frac{1}{\omega^j})}{1-\omega^j}\label{eq:firstderivpgfone}
\end{eqnarray}
Lastly, recall that $G_Y(1)=\sum_{k=0}^{\infty}p_k=1$. 

\subsubsection{Proof for \texorpdfstring{$\operatorname{Var}(U)$}{Var(U)}, Theorem \ref{theorem:expofU}}
From equation (\ref{eq:firstderivpgf}):
\begin{eqnarray}
n G^{''}_U(s) &=&  \frac{1}{4}G_{Y}(s)\sum^{n}_{t=1}(2r-4-n+2t)(2r-2-n+2t)s^{r-3-\frac{n}{2}+t}\nonumber\\
&&+ \frac{1}{2}G_{Y}^{'}(s)\sum^{n}_{t=1}s^{r-2-\frac{n}{2}+t} +\frac{1}{2}G_{Y}^{'}(s)\sum^{n}_{t=1}(2r-2-n+2t)s^{r-1-\frac{n}{2}+t}\nonumber\\ 
&& +\sum^{n}_{t=1}s^{r-1-\frac{n}{2}+t} G_{Y}^{''}(s)+\frac{1}{4}(n+2r-2)(n+2r)s^{\frac{n}{2}+r-1}\nonumber\\
&&-\frac{1}{4}(-n+2r-2)(-n+2r)s^{-\frac{n}{2}+r-1} \sum^{n-1}_{j=1}a(j)\frac{G_{Y}(\frac{s}{\omega^j})}{s-\omega^j}\nonumber\\
&& +\frac{1}{2}\resizebox{.9\hsize}{!}{$\Bigl((n+2r)s^{\frac{n}{2}+r-1}-(-n+2r)s^{-\frac{n}{2}+r-1}\Bigr) \sum^{n-1}_{j=1}a(j)\left(\frac{G_{Y}^{'}(\frac{s}{\omega^j})}{\omega^j(s-\omega^j)}-\frac{G_{Y}(\frac{s}{\omega^j})}{(s-\omega^j)^2}\right)$}\nonumber\\
&& +\frac{1}{2}\resizebox{.9\hsize}{!}{$\Bigl((n+2r)s^{\frac{n+2r}{2}}-(-n+2r)s^{\frac{-n+2r}{2}}\Bigr) \sum^{n-1}_{j=1}a(j)\left(\frac{G_{Y}^{'}(\frac{s}{\omega^j})}{\omega^j(s-\omega^j)}-\frac{G_{Y}(\frac{s}{\omega^j})}{(s-\omega^j)^2}\right)$}\nonumber\\
&& +\left(s^{\frac{n+2r}{2}}-s^{\frac{-n+2r}{2}}\right) \sum^{n-1}_{j=1}a(j)\left( \frac{G_{Y}^{''}(\frac{s}{\omega^j})}{\omega^{2j}(s-\omega^j)} -
2\frac{G_{Y}^{'}(\frac{s}{\omega^j})}{\omega^j(s-\omega^j)^2} + 2\frac{G_{Y}(\frac{s}{\omega^j})}{\left(s-\omega^j\right)^3}\right)\nonumber
\end{eqnarray}

Which leads to:
\begin{eqnarray}
n G^{''}_U(1) &=& \frac{1}{4}G_{Y}(1)\sum^{n}_{t=1}(2r-4-n+2t)(2r-2-n+2t)\nonumber\\
&&+G_{Y}^{'}(1)\sum^{n}_{t=1}(2r-2-n+2t)+n G_{Y}^{''}(1)\nonumber\\
&& +\frac{1}{4}(4n)\left(2r-1\right) \sum^{n-1}_{j=1}a(j)\frac{G_{Y}(\frac{1}{\omega^j})}{1-\omega^j}+ 2n \sum^{n-1}_{j=1}a(j)\left(\frac{G_{Y}^{'}(\frac{1}{\omega^j})}{\omega^j(1-\omega^j)} -\frac{G_{Y}(\frac{1}{\omega^j})}{\left(1-\omega^j\right)^2}\right)\nonumber\\
&=& (r^2+\frac{n^2}{4}-rn-3r+\frac{3n}{2}+2)n+(4r-2n-6)\frac{n(n+1)}{4}\nonumber\\
&&+\frac{n(n+1)(2n+1)}{6}G_{Y}(1)+(2r-1)n G_{Y}^{'}(1)+n G_{Y}^{''}(1)\nonumber\\
&& +n\left(2r-1\right) \sum^{n-1}_{j=1}a(j)\frac{G_{Y}(\frac{1}{\omega^j})}{1-\omega^j}+ 2n \sum^{n-1}_{j=1}a(j)\left(\frac{G_{Y}^{'}(\frac{1}{\omega^j})}{\omega^j(1-\omega^j)} -\frac{G_{Y}(\frac{1}{\omega^j})}{\left(1-\omega^j\right)^2}\right)\nonumber\\
&=& \frac{n}{12}\left(n^2+12r^2-24r+8\right)G_{Y}(1)+(2r-1)n G_{Y}^{'}(1)+n G_{Y}^{''}(1)+\nonumber\\
&& n\left(2r-1\right) \sum^{n-1}_{j=1}a(j)\frac{G_{Y}(\frac{1}{\omega^j})}{1-\omega^j}+ 2n \sum^{n-1}_{j=1}a(j)\left(\frac{G_{Y}^{'}(\frac{1}{\omega^j})}{\omega^j(1-\omega^j)} -\frac{G_{Y}(\frac{1}{\omega^j})}{\left(1-\omega^j\right)^2}\right)\label{eq:secondderivpgfone}
\end{eqnarray}
and 
\begin{eqnarray}
G^{''}_U(1) &=& \frac{1}{12}\left(n^2+12r^2-24r+8\right)+(2r-1)\operatorname{E}(Y)+ G_{Y}^{''}(1)+\nonumber\\
&& \left(2r-1\right) \sum^{n-1}_{j=1}a(j)\frac{G_{Y}(\frac{1}{\omega^j})}{1-\omega^j}+ 2 \sum^{n-1}_{j=1}a(j)\left(\frac{G_{Y}^{'}(\frac{1}{\omega^j})}{\omega^j(1-\omega^j)} -\frac{G_{Y}(\frac{1}{\omega^j})}{\left(1-\omega^j\right)^2}\right)\nonumber
\end{eqnarray}
Therefore $\operatorname{Var}(U)= G''_U(1)+\operatorname{E}(U)-(\operatorname{E}(U))^2$ becomes

\begin{eqnarray}
\operatorname{Var}(U) 
&=& G_{Y}^{''}(1)+\operatorname{E}(Y) - (\operatorname{E}(Y))^2+\frac{1}{12}(n^2-1)  -\left(2\operatorname{E}(Y)-1\right)\sum^{n-1}_{j=1}a(j)\frac{G_{Y}(\frac{1}{\omega^j})}{1-\omega^j}\nonumber\\
&& -\left(\sum^{n-1}_{j=1}a(j)\frac{G_{Y}(\frac{1}{\omega^j})}{1-\omega^j}\right)^2+2\sum^{n-1}_{j=1}a(j)\left(\frac{G_{Y}^{'}(\frac{1}{\omega^j})}{\omega^j(1-\omega^j)} -\frac{G_{Y}(\frac{1}{\omega^j})}{\left(1-\omega^j\right)^2}\right)\nonumber\\
&=& \operatorname{Var}(Y)+\frac{1}{12}(n^2-1)  -\left(2\operatorname{E}(Y)-1\right)\sum^{n-1}_{j=1}a(j)\frac{G_{Y}(\frac{1}{\omega^j})}{1-\omega^j}\nonumber\\
&& -\left(\sum^{n-1}_{j=1}a(j)\frac{G_{Y}(\frac{1}{\omega^j})}{1-\omega^j}\right)^2+2\sum^{n-1}_{j=1}a(j)\left(\frac{G_{Y}^{'}(\frac{1}{\omega^j})}{\omega^j(1-\omega^j)} -\frac{G_{Y}(\frac{1}{\omega^j})}{\left(1-\omega^j\right)^2}\right)\nonumber
\end{eqnarray}

\subsection{Proof of Corollary \ref{theorem:expupoissonN}}\label{sec:euvarpois}
For $\operatorname{E}(U)$ simply replace in (\ref{eq:expuPoi}) $G_{Y}(1/\omega^{j})$ and $\operatorname{E}(Y)$ by their respective values when $Y \sim Poisson(\theta)$.

\vspace{20pt}
For $\operatorname{Var}(U)$ replace in (\ref{eq:varuPoi}) $G_{Y}(1/\omega^{j}), G_{Y}^{'}(1/\omega^{j}), \operatorname{E}(Y)$, and $\operatorname{Var}(Y)$ by their respective values when $Y \sim Poisson(\theta)$ leads to
\begin{eqnarray}
	\operatorname{Var}(U)
	&=&\theta+\frac{1}{12}(n^2-1)-e^{-2\theta}\left(\sum^{n-1}_{j=1}a(j)\frac{e^{\frac{\theta}{\omega^j}}}{1-\omega^j}\right)^2\nonumber\\
	&&-e^{-\theta}(2\theta-1)\sum^{n-1}_{j=1}a(j)\frac{e^{\frac{\theta}{\omega^j}}}{1-\omega^j}\nonumber\\
	&&+2e^{-\theta}\sum^{n-1}_{j=1}a(j)\frac{e^{\frac{\theta}{\omega^j}}}{1-\omega^j}\left(\theta -\frac{1}{1-\omega^j}\right)\nonumber
	\end{eqnarray}

Re-expressing algebraically the answer is obtained.

\subsection{Proof of Lemma \ref{poissonlarge}}
\begin{proof}
In terms of the probabilities $p_k$, according to (\ref{eq:ln}) and (\ref{eq:ln2}) $\operatorname{E}(U)$ is given by,
\begin{equation}
   \operatorname{E}(U)= \sum_{l=0}^\infty\sum_{q=0}^{n-1}(n+ln)p_{n/2+a+ln+q}\nonumber
\end{equation}
where $a=0$ if $n$ is even and $a=1/2$ is $n$ is odd. At large values of $\lambda$, the Poisson probabilities are approximated by a Gaussian distribution as, 
\begin{equation}
     \operatorname{E}(U)\to \sum_{l=0}^\infty\sum_{q=0}^{n-1}(n+ln)\frac{e^{-(n/2+a+ln+q-\theta)^2/2\theta}}{\sqrt{2\pi\theta}}\nonumber
\end{equation}
We define the new variable $x=n/2+a+ln+q$, such that
\begin{equation}
    \operatorname{E}(U)\to J_1-J_2-J_3-J_4\nonumber
\end{equation}
where,
\begin{eqnarray}
J_1&\equiv& \sum_{x=0}^\infty x\,\frac{e^{-(x-\theta)^2/2\theta}}{\sqrt{2\pi\theta}},\nonumber\\
J_2&\equiv& \sum_{x=0}^{n/2+a} x\,\frac{e^{-(x-\theta)^2/2\theta}}{\sqrt{2\pi\theta}},\nonumber\\
J_3&\equiv& \sum_{l=0}^\infty\sum_{q=0}^{n-1}(a-n/2)\frac{e^{-(n/2+a+ln+q-\theta)^2/2\theta}}{\sqrt{2\pi\theta}},\nonumber\\
J_4&\equiv&\sum_{l=0}^\infty\sum_{q=0}^{n-1}q\frac{e^{-(n/2+a+ln+q-\theta)^2/2\theta}}{\sqrt{2\pi\theta}}\nonumber
\end{eqnarray}

The first sum is a standard Poisson-expected value, $J_1=\theta$.

For the second sum we can place an upper bound by substituting each term in the sum with the highest value of $x=n/2+a$,
\begin{equation}
J_2<(n/2+a)^2\frac{e^{-(n/2+a-\theta)^2/2\theta}}{\sqrt{2\pi\theta}}\nonumber
\end{equation}
which vanishes exponentially if $\theta\gg n/2+a$. The third integral is simply,
\begin{equation}
    J_3=(a-n/2)\sum_{x=0}^\infty \frac{e^{-(x-\theta)^2/2\theta}}{\sqrt{2\pi\theta}}=(a-n/2)\nonumber
\end{equation}
For the fourth integral, we can again place an upper bound by replacing each term in the sum over $q$ with its highest value,
\begin{equation}
    J_4< (n-1)\sum_{l=0}^\infty\sum_{q=0}^{n-1}\frac{e^{-(n/2+a+ln+q-\theta)^2/2\theta}}{\sqrt{2\pi\theta}}=(n-1)\nonumber
\end{equation}
Then it is easy to show,
\begin{equation}
    \lim_{\lambda\to\infty}\frac{1}{n\lambda}\operatorname{E}(U)=\lim_{\lambda\to\infty}\frac{1}{n\lambda}(J_1+J_2+J_3+J_4)=\lim_{\lambda\to\infty}\frac{1}{n\lambda}J_1=1\nonumber
\end{equation}

Now we turn to the variance, which can be written in the large $\lambda$ limit as 
\begin{eqnarray}
    \operatorname{Var}(U)&=&\sum_{l=0}^\infty\sum_{q=0}^{n-1}(n+ln)^2\frac{e^{-(n/2+a+ln+q-\theta)^2/2\theta}}{\sqrt{2\pi\theta}}\nonumber\\
    &&-\left(\sum_{l=0}^\infty\sum_{q=0}^{n-1}(n+ln)\frac{e^{-(n/2+a+ln+q-\theta)^2/2\theta}}{\sqrt{2\pi\theta}}\right)^2\nonumber
\end{eqnarray}
Again we define $x=n/2+a+ln+q$. such that
\begin{equation}
    \operatorname{Var}(U)=J_1-J_2-J_3-J_4+J_5-J_6+J_7\nonumber
\end{equation}
where,
\begin{eqnarray}
J_1&\equiv&\sum_{x=0}^\infty x^2\,\frac{e^{-(x-\theta)^2/2\theta}}{\sqrt{2\pi\theta}}-\left(\sum_{x=0}^\infty x\, \frac{e^{-(x-\theta)^2/2\theta}}{\sqrt{2\pi\theta}}\right)^2\nonumber\\
J_2&\equiv&\sum_{x=0}^{n/2+a} x^2\,\frac{e^{-(x-\theta)^2/2\theta}}{\sqrt{2\pi\theta}}-\left(\sum_{x=0}^{n/2+a} x\, \frac{e^{-(x-\theta)^2/2\theta}}{\sqrt{2\pi\theta}}\right)^2\nonumber\\
J_3&\equiv& \sum_{l=0}^\infty\sum_{q=0}^{n-1}2x(n/2-a-q)\frac{e^{-(n/2+a+ln+q-\theta)^2/2\theta}}{\sqrt{2\pi\theta}},\nonumber\\
J_4&\equiv& \sum_{l=0}^\infty\sum_{q=0}^{n-1}(n/2-a-q)^2\frac{e^{-(n/2+a+ln+q-\theta)^2/2\theta}}{\sqrt{2\pi\theta}},\nonumber\\
J_5&\equiv& \left(\sum_0^\infty 2x\frac{e^{-(x-\theta)^2/2\theta}}{\sqrt{2\pi\theta}}\right)\left(\sum_{l=0}^\infty\sum_{q=0}^{n-1}(n/2-a-q)\frac{e^{-(n/2+a+ln+q-\theta)^2/2\theta}}{\sqrt{2\pi\theta}}\right),\nonumber\\
J_6&\equiv& \left(\sum_0^{x=n/2+a} 2x\frac{e^{-(x-\theta)^2/2\theta}}{\sqrt{2\pi\theta}}\right)\left(\sum_{l=0}^\infty\sum_{q=0}^{n-1}(n/2-a-q)\frac{e^{-(n/2+a+ln+q-\theta)^2/2\theta}}{\sqrt{2\pi\theta}}\right),\nonumber\\
J_7&\equiv&\left(\sum_{l=0}^\infty\sum_{q=0}^{n-1}(n/2-a-q)\frac{e^{-(n/2+a+ln+q-\theta)^2/2\theta}}{\sqrt{2\pi\theta}}\right)^2\nonumber
\end{eqnarray}

The first expression is the standard Poisson variance, $J_1=\theta$. We again find bounds for the rest of the sums,
\begin{eqnarray}
J_2&<&(n/2+a)^3\frac{e^{-(n/2+a-\theta)^2/2\theta}}{\sqrt{2\pi\theta}}-\left((n/2+a)^2\frac{e^{-(n/2+a-\theta)^2/2\theta}}{\sqrt{2\pi\theta}}\right)^2\nonumber\\
J_4&<&(n/2-a-n+1)^2\nonumber\\
J_6&<&2(n/2-a+1)(n/2-a)^2(n-1)\frac{e^{-(n/2-a-\theta)^2/2\theta}}{\sqrt{2\pi\theta}}\nonumber\\
J_7&<&(n/2-a-n+1)^2\nonumber
\end{eqnarray}
and,
\begin{equation}
J_3-J_5=\sum_{l=0}^\infty\sum_{q=0}^{n-1}(2\theta-n/2+a+ln+q)q\frac{e^{-(n/2+a+ln+q-\theta)^2/2\theta}}{\sqrt{2\pi\theta}}.\label{threefive}
\end{equation}
In expression (\ref{threefive}) we define the new variable $y=ln/\theta$, and for large $\lambda$ we approximate the sum over $l$ with an integral over $y$, 
\begin{equation}
    J_3-J_5\approx \frac{\theta}{n}\int_0^\infty \sum_{q=0}^{n-1}2\theta(1-n/2\theta-a/\theta-y-q/\theta)q\,\frac{e^{-\theta(n/2\theta+a/\theta+y+q/\theta-1)^2/2}}{\sqrt{2\pi\theta}}dy\nonumber
\end{equation}
At large $\lambda$, the integral over $y$ can be calculated with a saddle point approximation, where the saddle point is given by $y_{sp}=1-n/2\theta-a/\theta-q/\theta$, resulting in,
\begin{equation}
    J_3-J_5\to 0+\mathcal{O}\left(\frac{1}{\sqrt{\theta}}\right)\nonumber
\end{equation}
Putting all these results together, we then find,
\begin{equation}
    \lim_{\lambda\to\infty}\frac{1}{n\lambda}\operatorname{Var}(U)=\lim_{\lambda\to\infty}\frac{1}{n\lambda}J_1=1\nonumber
\end{equation}
\end{proof}

\subsection{Proof for Result \ref{theorem:mlealpha}}
\begin{proof}
From (\ref{eq:poispmf}) we see that for $n>1$ the log-likelihood is,
\begin{eqnarray}
l &=& \log(L(\theta\mid u)) = (h(u)+g(u))\log(\theta)-\theta + \log(\sum_{q=0}^{n-1-g(u)}\frac{\theta^{q}}{(h(u) + g(u)+q)!})\nonumber\\
&=& c\log(\theta) - \theta +
\log(\sum_{q=0}^{n-1-g(u)}d_q\theta^q)\nonumber
\end{eqnarray}

where $c=h(u)+g(u)$ and $d_q=\frac{1}{(c+q)!}$. Then
\begin{eqnarray}
\frac{dl}{d \theta}&=&\frac{c}{\theta} -1+\frac{\sum_{q=1}^{n-1-g(u)}d_qq\theta^{q-1}}{\sum_{q=0}^{n-1-g(u)}d_q\theta^q}=0\nonumber\\
&=&\frac{c-\theta}{\theta}+\frac{\sum_{q=1}^{n-1-g(u)}d_qq\theta^{q-1}}{\sum_{q=0}^{n-1-g(u)}d_q\theta^q}=0\nonumber\\
&=&\frac{c\sum_{q=0}^{n-1-g(u)}d_q\theta^q -\theta\sum_{q=0}^{n-1-g(u)}d_q\theta^q+\theta\sum_{q=1}^{n-1-g(u)}d_qq\theta^{q-1}}{\theta\sum_{q=0}^{n-1-g(u)}d_q\theta^q}=0\nonumber\\
&=&c\sum_{q=0}^{n-1-g(u)}d_q\theta^q -\theta\sum_{q=0}^{n-1-g(u)}d_q\theta^q+\theta\sum_{q=1}^{n-1-g(u)}d_qq\theta^{q-1}=0\nonumber\\
&=&c\sum_{q=0}^{n-1-g(u)}d_q\theta^q -\sum_{q=0}^{n-1-g(u)}d_q\theta^{q+1}+\sum_{q=1}^{n-1-g(u)}d_qq\theta^{q}=0\nonumber\\
&=&c\sum_{q=0}^{n-1-g(u)}d_q\theta^q -\sum_{q=1}^{n-1+1-g(u)}d_{q-1}\theta^{q}+\sum_{q=1}^{n-1-g(u)}d_qq\theta^{q}=0\nonumber\\ 
&=&cd_0-d_{n-1}\theta^{n-1+1-g(u)}+c\sum_{q=1}^{n-1-g(u)}d_{q}\theta^q -\sum_{q=1}^{n-1-g(u)}(c+q)d_{q}\theta^{q}+\nonumber\\
&&\sum_{q=1}^{n-1-g(u)}d_qq\theta^{q}=0 \nonumber\\
&=&cd_0-d_{n-1}\theta^{n-1+1}=0\nonumber 
\end{eqnarray}

The equality before last occurs because $d_{q-1} = (c+q)d_{q}$. Therefore,
\begin{eqnarray}
\hat{\theta}&=&(\frac{cd_0}{d_{n-1}})^{1/(n-1+1)}=\Bigl((h(u)+g(u)+n-1)....(h(u)+g(u))\Bigr)^{1/n}\nonumber\\
&=&\begin{cases}
\Bigl((u-\frac{n}{2}+\frac{1}{2}+g(u)+n-1)...(u-\frac{n}{2}+\frac{1}{2}+g(u))\Bigr)^{\frac{1}{n}} \hspace{0.5cm} n \text{ is odd}\nonumber\\
\Bigl((u+\frac{n}{2}-1+g(u))....(u-\frac{n}{2}+g(u))\Bigr)^{\frac{1}{n}} \hspace{0.5cm} n \text{ is even}
\end{cases}\nonumber\\
&=&\begin{cases}
\prod_{q=0}^{n-1-g(u)} \Bigl(u-\frac{n}{2}+\frac{1}{2}+g(u)+q\Bigr)^{\frac{1}{n}} \hspace{0.5cm} n \text{ is odd}\nonumber\\
\prod_{q=0}^{n-1-g(u)} \Bigl(u-\frac{n}{2}+g(u)+q\Bigr)^{\frac{1}{n}} \hspace{0.5cm} n \text{ is even}\nonumber\\
\end{cases}\\
&=&\prod_{q=0}^{n-1-g(u)} \Bigl(\lceil u -\frac{n}{2}\rceil+g(u)+q\Bigr)^{\frac{1}{n}}\nonumber
\end{eqnarray}
This MLE adjusts for the effect of rounding to the nearest integer. Occasions when $\lceil u -\frac{n}{2}\rceil+g(u)+q<0$ are of probability zero and thus must be omitted before calculating the geometric mean,

\begin{eqnarray}
\hat{\theta}
&=&\prod_{q \in \mathcal{P}} \Bigl(\lceil u -\frac{n}{2}\rceil+g(u)+q\Bigr)^{\frac{1}{m}}\nonumber
\end{eqnarray}

where $\mathcal{P}$ is the set such that $\lceil u -\frac{n}{2}\rceil+g(u)+q>0$ and $m$ is the length of $\mathcal{P}$.
\end{proof}

\subsection{Proof of Result \ref{eq:PoissonMLEthetalarge}}
\begin{proof}

Considering large $\lambda$ we omit $g(u)$ from (\ref{eq:mlealpha}),
\begin{equation}
\operatorname{E}(\hat\theta)=\sum_{u/n=0}^{\infty} \prod_{k\in\mathcal{P}}\Bigl(\lceil u-n/2\rceil+k\Bigr)^{\frac{1}{m}}P(U=u)\nonumber
\end{equation}
At large $\lambda$, with (\ref{eq:ln}) and (\ref{eq:ln2}) we use the fact that $P(Y=y)$ is approximately Gaussian, and we can express the expectation value as,
\begin{equation}
\operatorname{E}(\hat\theta)=\sum_{l=0}^\infty\sum_{q=0}^{n-1}\prod_{k\in\mathcal{P}}(ln-n/2+a)\prod_{k\in\mathcal{P}}\left(1+\frac{k}{ln-n/2+a}\right)^{\frac{1}{m}}\frac{e^{-(ln+n/2+a+q-\theta)^2/2\theta}}{\sqrt{2\pi\theta}}\nonumber
\end{equation}

where $a=0$ if $n$ is even and $a=1/2$ is $n$ is odd. The Gaussian distribution for large $\lambda$ implies that the only contributions to the sum that are not exponentially small are where $ln\sim\theta$, therefore we can assume $ln\gg n, a, k, q$, in the expression for $\operatorname{E}(\hat\theta)$ and,
\begin{equation}
\operatorname{E}(\hat\theta)\approx \sum_{l=1}^\infty\sum_{q=0}^{n-1}ln\prod_{k=0}^{n-1} \left(1+\frac{k}{ln}\right)^{\frac{1}{n-1}}\frac{e^{-(ln+n/2+a+q-\theta)^2/2\theta}}{\sqrt{2\pi\theta}}\nonumber
\end{equation}
where since we exclude the probability that $l=0$, now all $\mathcal{P}$ are of the same size $n-1$. Furthermore, since $k\ll ln$, we can approximate,
\begin{eqnarray}
\operatorname{E}(\hat\theta)&\approx& \sum_{l=0}^\infty\sum_{q=0}^{n-1}ln\frac{e^{-(ln+n/2+a+q-\theta)^2/2\theta}}{\sqrt{2\pi\theta}}\nonumber\\
&=&\operatorname{E}(U)-\sum_{l=0}^\infty\sum_{q=0}^{n-1} n\frac{e^{-(ln+n/2+a+q-\theta)^2/2\theta}}{\sqrt{2\pi\theta}}\nonumber
\end{eqnarray}
Now invoking Lemma \ref{poissonlarge}, we find 
\begin{equation}
\lim_{\lambda\to\infty}\frac{1}{n\lambda}\operatorname{E}(\hat\theta)=\lim_{\lambda\to\infty}\frac{1}{n\lambda}\left(\operatorname{E}(U)-n\right)=1\nonumber
\end{equation}
\end{proof}

\subsection{Proof of Result \ref{eq:PoissonMLEnlarge}}
\begin{proof}

The expected value of the MLE is given by,
\begin{equation}
\operatorname{E}(\hat\theta)=\sum_{u/n=0}^{\infty} \prod_{k\in\mathcal{P}}\left(\lceil u-n/2\rceil+k\right)^{\frac{1}{m}}P(U=u)\nonumber
\end{equation}
where we express the sum in terms of $u/n$ which has a non-negative integers support. From Lemma \ref{theorem:pmfu}, for large $n$ we have $v_0 = \lfloor \lambda+0.5\rfloor$ and $P(U=u) \approx \delta^{\rm Kronecker}_{u,nv_0}$, so  
\begin{equation}
\lim_{n\to\infty}\operatorname{E}(\hat\theta)=\lim_{n\to\infty} \prod_{k=0}^{n-1}(\lceil nv_0-n/2\rceil+k)^{\frac{1}{n}}=\lim_{n\to\infty}\lceil nv_0-n/2\rceil\prod_{k=0}^{n-1}\left(1+\frac{k}{\lceil nv_0-n/2\rceil}\right)^{\frac{1}{n}}\nonumber
\end{equation}

In the large $n$ limit, It is more convenient to work with the logarithm of $\operatorname{E}(\hat\theta)$,
\begin{equation}
\lim_{n\to\infty}\log(\operatorname{E}(\hat\theta))=\lim_{n\to\infty}\left[\log(\lceil nv_0-n/2\rceil)+\sum_{k=0}^{n-1}\frac{1}{n}\log\left(1+\frac{k}{\lceil nv_0-n/2\rceil}\right)\right]\nonumber
\end{equation}
In the large $n$ limit, the sum over $k$ can be approximated with an integral over variable $x=k/\lceil nv_0-n/2\rceil$, 
\begin{eqnarray}
\lim_{n\to\infty}\log(\operatorname{E}(\hat\theta))&=&\lim_{n\to\infty}\left[\log(\lceil nv_0-n/2\rceil)+\frac{\lceil nv_0-n/2\rceil}{n}\int_0^\frac{n-1}{\lceil nv_0-n/2\rceil}\log(1+x)dx\right]\nonumber\\
&=& \resizebox{.70\hsize}{!}{$\lim_{n\to\infty}\left[\log(\lceil nv_0-n/2\rceil)+\Bigl(v_0-\frac{1}{2}\Bigr)\Bigl(\frac{1}{v_0-\frac{1}{2}}+1\Bigr)\log\left(\frac{1}{v_0-1/2}+1\right)-1\right]$}\nonumber
\end{eqnarray}
Exponentiating both sides and dividing by $n$ we get our result (\ref{largenmse}).

\end{proof}

\subsection{Proof of Result \ref{eq:PoissonMLEnlargev2}}
\begin{proof}

Observe that in the large $n$ limit, probability mass concentrates more in $u=0$, such that 
\begin{eqnarray}
\lim_{n\to\infty}\operatorname{E}(\hat\theta)=\lim_{n\to\infty}\prod_{k=n/2}^{n-1}(\lceil-n/2\rceil+k)^{2/n}=\lim_{n\to\infty}\frac{n}{2}\prod_{k=n/2}^{n-1}\left(-1+\frac{2k}{n}\right)^{2/n}\nonumber
\end{eqnarray}

We now take the logarithm on both sides, such that
\begin{eqnarray}
\lim_{n\to\infty}\log(\operatorname{E}(\hat\theta))=\lim_{n\to\infty}\left[\log\left(\frac{n}{2}\right)+\sum_{k=n/2}^{n-1}\frac{2}{n}\log\left(-1+\frac{2k}{n}\right)\right]\nonumber
\end{eqnarray}
where the sum can be approximated by an integral in the large $n$ limit as,
\begin{eqnarray}
\lim_{n\to\infty}\log(\operatorname{E}(\hat\theta))=\lim_{n\to\infty}\left[\log\left(\frac{n}{2}\right)+\int_1^2 \log\left(-1+x\right)dx\right]=\log\left(\frac{n}{2}\right)-1\nonumber
\end{eqnarray}
Exponentiating both sides of the equation, and dividing by $n$, we have,
\begin{equation}
\lim_{n\to\infty}\frac{1}{n}\operatorname{E}(\hat\theta)=\frac{1}{2e}\nonumber
\end{equation}
\end{proof}

\subsection{Proof for Corollary \ref{theorem:binomU}}

For $\operatorname{E}(U)$ simply replace $G_{Y}(1/\omega^{j})$ and $\operatorname{E}(Y)$ by their respective values when 

$Y \sim binomial(mn,\phi)$.

For $\operatorname{Var}(U)$:
\begin{equation}
\begin{split}
n G''_U(s) =& \frac{mn\phi}{2}\left(q+\phi s\right)^{mn-1}\sum^{n}_{t=1}(2r-n-2+2t)s^{r-2+t-\frac{n}{2}}+\nonumber\\
& \frac{1}{2}\left((q+\phi s)^mn-q^mn\right)\sum^{n}_{t=1}
\left(r-2+t-\frac{n}{2}\right)\left(2r-n-2+2t\right)s^{r-3+t-\frac{n}{2}}+\nonumber\\
& mn(mn-1)\phi^2\left(q+\phi s\right)^{mn-2}\sum^{n}_{t=1}
s^{r-1+t-\frac{n}{2}}+\nonumber\\
& mn\phi\left(q+\phi s\right)^{mn-1}\sum^{n}_{t=1}\left(r-1+t-\frac{n}{2}\right)s^{r-2+t-\frac{n}{2}}\nonumber\\ 
& +\frac{1}{2}(\left(\frac{n}{2}+r-1\right)\left(n+2r\right)s^{\frac{n}{2}+r-2}\nonumber\\
& - \left(-\frac{n}{2}+r-1\right)\left(-n+2r\right)s^{-\frac{n}{2}+r-2}) \sum^{n-1}_{j=1}\frac{(q+\frac{\phi s}{\omega^j})^{mn}-q^{mn}}{s-\omega^j}\nonumber\\
& \resizebox{.9\hsize}{!}{$+\left((n+2r)s^{\frac{n+2r-2}{2}}-(-n+2r)s^{\frac{-n+2r-2}{2}}\right) \sum^{n-1}_{j=1}\left(\frac{\frac{mn\phi}{\omega^j}(q+\frac{\phi s}{\omega^j})^{mn-1}}{s-\omega^j} -\frac{(q+\frac{\phi s}{\omega^j})^{mn}-q^{mn}}{\left(s-\omega^j\right)^2}\right)$}\nonumber\\
& \resizebox{.9\hsize}{!}{$+\left(s^{\frac{n+2r}{2}}-s^{\frac{-n+2r}{2}}\right) \sum^{n-1}_{j=1}(\frac{mn(mn-1)\left(\frac{\phi}{\omega^j}\right)^2(q+\frac{\phi s}{\omega^j})^{mn-2}}{s-\omega^j}-2\frac{\frac{mn\phi}{\omega^j}(q+\frac{\phi s}{\omega^j})^{mn-1}}{\left(s-\omega^j\right)^2}$}\nonumber\\
&+2\frac{(q+\frac{\phi s}{\omega^j})^{mn}-q^{mn}}{\left(s-\omega^j\right)^3})\nonumber
\end{split}
\end{equation}

Thus, $n G''_U(1)$ is given by:
\begin{eqnarray} 
n G''_U(1) &=& \frac{mn\phi}{2}n\left(2r-1\right)+ \frac{n}{12}\left(n^2+12r^2-24r+8\right)\left(1-q^{mn}\right) \nonumber\\
&&+mn(mn-1)n\phi^2+(r-\frac{1}{2})m n^2\phi\nonumber\\
&& -n\left(2r-1\right) \sum^{n-1}_{j=1}\frac{(q+\frac{\phi}{\omega^j})^{mn}-q^{mn}}{\omega^j-1}\nonumber\\
&&- 2n \sum^{n-1}_{j=1}\left(\frac{\frac{mn\phi}{\omega^j}(q+\frac{\phi}{\omega^j})^{mn-1}}{\omega^j-1} +\frac{(q+\frac{\phi}{\omega^j})^{mn}-q^{mn}}{\left(\omega^j-1\right)^2}\right)\nonumber
\end{eqnarray}

Some additional algebraic calculations lead to the result.
\end{appendices}

 \bibliography{bibtexreferences}

\end{document}